\documentclass[11pt]{article}
\usepackage[margin=1.05in]{geometry}
\usepackage[T1]{fontenc}
\usepackage{lmodern,microtype}
\usepackage{amsmath,amssymb,amsthm,mathtools,bm,mathrsfs}
\usepackage{enumitem,booktabs,array}
\usepackage{xcolor,hyperref}
\hypersetup{colorlinks=true,linkcolor=blue!45!black,citecolor=blue!45!black,urlcolor=blue!45!black,
  pdftitle={Angular Besov Refinement of the Spherical Stein--Tomas Theorem},
  pdfauthor={Rishad Shahmurov and Veli Shahmurov},
  pdfkeywords={Fourier restriction, angular Besov spaces, nonlinear Helmholtz equation, scattering, spherical harmonics}}
\newtheorem{theorem}{Theorem}[section]

\newtheorem{corollary}[theorem]{Corollary}
\newtheorem{lemma}[theorem]{Lemma}
\newtheorem{remark}[theorem]{Remark}

\numberwithin{equation}{section}
\newcommand{\Sph}{\mathbb S^2}
\newcommand{\RR}{\mathbb R}
\newcommand{\HH}{\mathcal H}
\newcommand{\LL}{\mathbf L}

\newcommand{\cN}{\mathcal N}
\newcommand{\supp}{\operatorname{supp}}
\newcommand{\dist}{\operatorname{dist}}
\newcommand{\ang}{\operatorname{ang}}

\newcommand{\norm}[1]{\left\lVert#1\right\rVert}

\newcommand{\pair}[2]{\left\langle#1,#2\right\rangle}
\newcommand{\one}{\mathbf1}
\title{\bfseries Angular Besov Refinement of the Spherical Stein--Tomas Theorem}
\author{Rishad Shahmurov\thanks{Cellular Products Research and Development, 595 East Crossville Rd., Roswell, GA 30075, USA. Email: \texttt{rshahmurov@crimson.ua.edu}.}
\and Veli Shahmurov\thanks{Antalya Bilim University, Department of Industrial Engineering, 07190 Antalya, Turkey; Azerbaijan State Economic University, AZ1001 Baku, Azerbaijan. Emails: \texttt{veli.sahmurov@antalya.edu.tr}, \texttt{veli.sahmurov@gmail.com}.}}
\date{}
\begin{document}
\maketitle
\begin{abstract}
The spherical Stein--Tomas inequality maps $L^2(\mathbb S^2)$ into
$L^4(\mathbb R^3)$. We study its refinement in which the angular input
is measured by the strictly weaker Besov norm $B^0_{2,4}(\mathbb S^2)$.
The central separated-band inequality bounds the squared $L^2$ norm of
the product of extensions from degrees $L\ll K$ by $C(L/K)$ times
the squared input norms. An exact spherical-Hankel expansion controls
the far field, while angular-momentum leakage and transversality handle
the Fresnel region. We construct scale-separated, normalized angular
bands for which the extension norm is equivalent to the $\ell^4$ norm
of the coefficients. This proves that the fourth-power summation is
optimal and quantitatively distinguishes the result from the classical
$\ell^2$ summation. As an application, small angular Besov data generate
small $L^4$ solutions to the cubic nonlinear Helmholtz equation in
$\mathbb R^3$, with a smallness condition that remains sharp on multiscale data even
when their $L^2$ norm is arbitrarily large or infinite. The associated on-shell nonlinear scattering amplitude has
improved angular $B^0_{2,4/3}$ regularity, and its leading cubic term
is quantitatively decoupled across widely separated angular bands.
\end{abstract}
\medskip
\noindent\textbf{Keywords.} Fourier restriction; Stein--Tomas inequality; angular Besov space; sharp multiscale summation; nonlinear Helmholtz equation; scattering.\par
\medskip
\noindent\textbf{2020 Mathematics Subject Classification.} Primary 42B10, 42B37; Secondary 42C10, 35J60, 46E35.

\section{Introduction}
\subsection{Why a finer angular norm is needed}
The Fourier transform of a function on $\mathbb R^3$ is not ordinarily defined on a set of three-dimensional Lebesgue measure zero. Fourier restriction theory asks when one can nevertheless give meaning to its values on a curved surface. For the unit sphere, the dual formulation is the extension operator
\begin{equation}\label{eq:Edef}
 Eg(x)=\int_{\Sph}e^{ix\cdot\omega}g(\omega)\,d\sigma(\omega),\qquad x\in\RR^3,
\end{equation}
where $d\sigma$ is surface measure of total mass $4\pi$. The endpoint Stein--Tomas estimate in three dimensions is
\begin{equation}\label{eq:ST}
 \norm{Eg}_{L^4(\RR^3)}\le C_{\mathrm{ST}}\norm{g}_{L^2(\Sph)}.
\end{equation}
It is one of the basic curved-surface estimates in harmonic analysis: curvature creates enough cancellation for a superposition of plane waves to have a global $L^4$ bound.

Estimate \eqref{eq:ST} treats every angular harmonic on the same footing and aggregates their sizes in $\ell^2$. In a problem with many widely separated angular scales, this may overstate the amount of input regularity needed. The question we study is whether one may aggregate the contributions from \emph{dyadic angular frequencies} in the larger sequence space $\ell^4$ instead. This is an endpoint question: there is no room to improve the spatial exponent $4$ by an elementary Sobolev embedding, and a generic Littlewood--Paley decomposition alone does not supply the required cross-scale decay.

For concreteness, write $\mathcal H_k$ for the spherical harmonics of degree $k$. With $Q=-\Delta_{\Sph}$, set $\Lambda=(Q+\tfrac14)^{1/2}$, so that $\Lambda$ has eigenvalue $k+\tfrac12$ on $\mathcal H_k$. Choose $\varphi_0\in C_c^\infty([0,2))$, equal to one on $[0,1]$, and define the inhomogeneous angular decomposition
\begin{equation}\label{eq:LP}
 P_0=\varphi_0(\Lambda),\qquad
 \varphi(s)=\varphi_0(s)-\varphi_0(2s),\qquad
 P_j=\varphi(2^{-j}\Lambda),\quad j\ge1.
\end{equation}
The corresponding Besov norm is
\begin{equation}\label{eq:Besov}
 \norm{g}_{B^0_{2,4}(\Sph)}\asymp
 \left(\sum_{j\ge0}\norm{P_jg}_{L^2(\Sph)}^4\right)^{1/4}.
\end{equation}
Different smooth inhomogeneous decompositions yield equivalent norms. Unlike the $L^2$ norm, \eqref{eq:Besov} allows the total squared harmonic energy to diverge, provided its contributions on dyadic bands have a fourth-power summability property. Thus $L^2(\Sph)$ is continuously and strictly contained in $B^0_{2,4}(\Sph)$.

The gain matters most for data spread across many frequency scales.
In Section~\ref{sec:multiscale} we construct normalized functions $h_m$
with mutually separated angular bands such that, for every finite scalar
sequence $(c_m)$,
\begin{equation}\label{eq:intro-exact-law}
 c\left(\sum_m|c_m|^4\right)^{1/4}
 \le \norm{E\sum_m c_mh_m}_{L^4(\RR^3)}
 \le C\left(\sum_m|c_m|^4\right)^{1/4}.
\end{equation}
In contrast, the classical Stein--Tomas bound uses
$(\sum_m|c_m|^2)^{1/2}$. For $J$ equally weighted scales, this is a loss
of order $J^{1/4}$. In particular the enlarged Besov domain is not
merely formal: \eqref{eq:intro-exact-law} produces distributional
angular densities outside $L^2(\Sph)$ whose extensions have finite
$L^4$ norm, with the conclusion not obtainable from the triangle
inequality on the individual extensions.

The principal application concerns the stationary cubic nonlinear
Helmholtz equation
\begin{equation}\label{eq:intro-nonlinear}
 (-\Delta-1)u=Q(x)|u|^2u\qquad\text{in }\RR^3.
\end{equation}
The cubic power is exactly matched to the endpoint $L^4$: if
$u\in L^4$, then $|u|^2u\in L^{4/3}$. An outgoing Helmholtz
resolvent maps $L^{4/3}$ into $L^4$, while the refined extension
bound supplies incident fields controlled by a much larger angular
data class. We show that a nonlinear solution may be small in
$L^4$ although its incoming angular density has arbitrarily large,
or infinite, $L^2$ norm. The dual restriction estimate gives additional
information about the angular distribution of its scattered radiation.
The point is not that a large $L^2$ norm by itself prevents a small
Helmholtz wave: a single high-degree harmonic may already have a small
extension norm. Rather, on the constructed family of nearly saturating
angular bands the refined theorem identifies the \emph{correct}
multiscale smallness condition. With $J$ equal coefficients the incident
wave has size comparable to $J^{1/4}$, whereas using Stein--Tomas on
its total $L^2$ datum produces $J^{1/2}$. This quantitative distinction,
not the mere largeness of $L^2$, is the application below.

\subsection{Relation to earlier restriction and spectral theory}
Tomas's restriction theorem, Stein's endpoint argument, and Strichartz's work linking restriction to dispersive equations established the classical $L^2$ theory on curved surfaces \cite{tomas1975,stein1993,strichartz1977}. Foschi subsequently identified global extremizers for the sphere extension inequality \cite{foschi2015}. Our aim is different from sharpening the numerical Stein--Tomas constant: we retain the same spatial $L^4$ norm but enlarge the class of admissible angular data.

There are several distinct ways to refine restriction estimates. Bak and Seeger proved Lorentz-space refinements in the \emph{spatial} variable and for more general measures and oscillatory integrals \cite{bakseeger2011}. Here the refinement is instead in the \emph{angular spectral summation} of the input on the fixed sphere. Smooth spectral multipliers and Besov spaces on spheres were developed by Strichartz, Bonami--Clerc, and Narcowich--Petrushev--Ward \cite{strichartz1972,bonamiclerc1973,npwSIAM2006,npwJFA2006}; these supply the framework in which the new norm is intrinsic. The proof additionally uses a local eigenfunction $L^2$ estimate of Sogge \cite{sogge2016} and an Agmon--H\"ormander local extension estimate \cite{agmonhormander1976}.

Bilinear restriction theory emphasizes how separated oscillatory contributions interact; Foschi--Klainerman and Burq--G\'erard--Tzvetkov supply useful neighboring viewpoints \cite{foschiklainerman2000,bgt2005}. Our separated-band problem differs from the standard separated-cap problem. A complete spherical-harmonic band is not supported in a single small angular cap, and its physical extension moves through a Bessel turning-point region. We therefore use the cap estimate only \emph{after} proving a quantitative localization statement. Estimates for large spheres and folding canonical relations, such as those of Bennett--Seeger, concern related but different oscillatory regimes \cite{bennettseeger2009}.

A different line of work obtains restriction and Strichartz refinements
by testing how the input concentrates on small \emph{directional or
frequency caps}. Bourgain's restriction and $X_p$-type ideas for
nonlinear dispersive problems, the estimates of Moyua--Vargas--Vega,
and the refinements developed by B\'egout--Vargas are landmarks of
this viewpoint \cite{bourgain1993,mvv1999,begoutvargas2007}.
Those contributions concern related dispersive restriction geometries,
not the precise angular spectral decomposition of the present theorem.
A small cap singles out a direction of propagation, whereas a dyadic
spherical-harmonic band singles out the magnitude of angular momentum
about a spatial center. Neither type of localization can simply be
substituted for the other in the argument.

For the nonlinear application, Guti\'errez developed small-data $L^4$ solutions of the cubic Helmholtz problem from $L^2$ Herglotz data \cite{gutierrez2004}. Later work of Gell-Redman, Hassell, Shapiro and Zhang established asymptotic expansions for nonlinear Helmholtz eigenfunctions in higher-regularity scattering frameworks \cite{gellredman2020}. Symmetry-improved Fourier extension estimates have also been applied to nonlinear Helmholtz equations \cite{wethyesil2021}. Our use of a zero-order angular Besov input, rather than improved spatial integrability or symmetry, addresses a different aspect of the incident data.

For the linear Helmholtz discussion, there is also a substantial established theory of Herglotz waves and distributional spherical densities. In particular, Barcel\'o, Folch-Gabayet, Luque, P\'erez-Esteva and Vilela characterize entire Helmholtz solutions arising from Sobolev-distribution densities \cite{barcelo2021}. We do not claim that distributional Herglotz densities are new. The additional assertion here is the \emph{specific global $L^4$ criterion} furnished by the zero-order Besov norm, together with an explicit example outside the traditional $L^2$ Herglotz class.

\subsection{Main results and their meaning}
Our first theorem enlarges the natural domain of the spherical extension operator without changing its spatial endpoint.
\begin{theorem}\label{thm:main}
For every $g\in B^0_{2,4}(\Sph)$, the extension operator \eqref{eq:Edef}, initially defined on finite harmonic sums, has a unique extension to $L^4(\RR^3)$ and satisfies
\begin{equation}\label{eq:main}
 \boxed{\quad \norm{Eg}_{L^4(\RR^3)}\le C\norm{g}_{B^0_{2,4}(\Sph)}.\quad}
\end{equation}
The constant depends only on the smooth angular decomposition and not on the number of active harmonic bands.
\end{theorem}

The main new estimate behind this result isolates the mechanism governing interactions between very different angular degrees. It is stated for entire spectral bands, not for a single pair of spherical harmonics.
\begin{theorem}\label{thm:bilinear}
There exist constants $C,C_0>0$ such that if $L,K\ge1$, $K\ge C_0L$, and
\begin{equation}\label{eq:bands}
 \supp_{\Lambda}g_L\subset[L/2,2L],\qquad
 \supp_{\Lambda}g_K\subset[K/2,2K],
\end{equation}
then
\begin{equation}\label{eq:bilinear}
 \boxed{\quad\norm{Eg_LEg_K}_{L^2(\RR^3)}^2
 \le C\frac LK\norm{g_L}_2^2\norm{g_K}_2^2.\quad}
\end{equation}
The power $L/K$ is optimal in the sense of uniform bounds over separated bands.
\end{theorem}

Why does \eqref{eq:bilinear} prove a Besov estimate? At a fixed radius the extension preserves spherical-harmonic degree. The angular Littlewood--Paley square-function inequality converts the fourth power of the spatial norm into a double sum of products $Eg_jEg_k$. For neighboring bands, Stein--Tomas suffices. For distant bands, \eqref{eq:bilinear} produces a summable factor $2^{-|j-k|}$. Discrete convolution then replaces the original $\ell^2$ input summation by $\ell^4$. We give all steps in Section~\ref{sec:spectral}.

Duality supplies a genuinely stronger restriction target as well. We use the convention $\widehat f(\xi)=\int e^{-ix\cdot\xi}f(x)\,dx$.
\begin{corollary}\label{cor:dual}
For every $f\in L^{4/3}(\RR^3)$, its Fourier restriction to $\Sph$ belongs to $B^0_{2,4/3}(\Sph)$ and obeys
\begin{equation}\label{eq:dual}
 \norm{\widehat f|_{\Sph}}_{B^0_{2,4/3}(\Sph)}
 \le C\norm{f}_{L^{4/3}(\RR^3)}.
\end{equation}
\end{corollary}
\begin{proof}
We give the dyadic duality argument explicitly, including the low
angular frequencies. This is useful because, for a general function in
$L^{4/3}(\RR^3)$, the expression $\widehat f|_{\Sph}$ is initially a
\emph{restriction distribution}, not necessarily the pointwise trace
of a continuous Fourier transform.

\smallskip\noindent\emph{Orthogonal angular blocks.}
Let $D_{-1}$ be the orthogonal projector onto the eigenspaces of
$\Lambda$ with eigenvalues less than $1$, and, for $j\ge0$, set
\[
 D_j=\one_{[2^j,2^{j+1})}(\Lambda).
\]
The $D_j$ have finite-dimensional, mutually orthogonal ranges and
$\sum_{j\ge-1}D_j=I$ on $L^2(\Sph)$. The sharp blocks are introduced
only for this proof; the statement of the theorem continues to use the
smooth projections $P_j$ in \eqref{eq:LP}. For every $1\le q<\infty$,
\begin{equation}\label{eq:dual-sharp-equivalence}
 c_q\Big(\sum_{j\ge-1}\norm{D_jh}_2^q\Big)^{1/q}
 \le \norm h_{B^0_{2,q}(\Sph)}
 \le C_q\Big(\sum_{j\ge-1}\norm{D_jh}_2^q\Big)^{1/q}.
\end{equation}
Here and below the same formula holds for distributions for which the
right side is finite. To see \eqref{eq:dual-sharp-equivalence} directly,
observe that the spectral support of each $P_k$ meets only a fixed
number of the intervals defining the $D_j$, and conversely each $D_j$
meets only a fixed number of the supports of the $P_k$. Thus, by
orthogonality and the boundedness of the scalar multipliers,
\[
 \norm{P_kh}_2\le C\sum_{|j-k|\le C_1}\norm{D_jh}_2.
\]
For the reverse estimate use $\sum_{k\ge0}P_k=I$ on each finite
harmonic sum and write
\[
 D_jh=\sum_{|k-j|\le C_1}D_jP_kh,
 \qquad
 \norm{D_jh}_2\le\sum_{|k-j|\le C_1}\norm{P_kh}_2.
\]
There are only finitely many exceptional low indices, absorbed into
the constant $C_1$. Taking the $\ell^q$ norm in these two finite
convolutions proves \eqref{eq:dual-sharp-equivalence}. By harmonic
approximation the equivalence extends to the asserted distributional
spaces.

\smallskip\noindent\emph{The estimate for Schwartz data.}
Take $f\in\mathcal S(\RR^3)$, so that $h=\widehat f|_{\Sph}$ is a
smooth function. Let $g$ be any finite spherical-harmonic sum. Our
Fourier-transform convention and Fubini's theorem give
\begin{equation}\label{eq:dual-pairing-expanded}
 \begin{aligned}
 \int_{\Sph} h(\omega)\overline{g(\omega)}\,d\sigma(\omega)
 &=\int_{\RR^3} f(x)
   \overline{\int_{\Sph}e^{ix\cdot\omega}g(\omega)
              \,d\sigma(\omega)}\,dx\\
 &=\int_{\RR^3}f(x)\overline{Eg(x)}\,dx.
 \end{aligned}
\end{equation}
By H\"older's inequality and Theorem~\ref{thm:main},
\begin{equation}\label{eq:dual-pairing-estimate}
 \left|\int_{\Sph}h\overline g\,d\sigma\right|
 \le\norm f_{4/3}\norm{Eg}_4
 \le C\norm f_{4/3}\norm g_{B^0_{2,4}}
 \le C'\norm f_{4/3}
        \Big(\sum_{j\ge-1}\norm{D_jg}_2^4\Big)^{1/4}.
\end{equation}

We now identify the complete target norm, rather than invoking an
abstract duality statement. Fix a finite set $J$ of block indices and
put $b_j=\norm{D_jh}_2$. Whenever $b_j>0$, let
$e_j=D_jh/b_j$, which has norm one and lies in the range of $D_j$;
for $b_j=0$ omit the corresponding index. Given nonnegative numbers
$a_j$, define $g=\sum_{j\in J}a_je_j$. Orthogonality gives
\[
 \norm{D_jg}_2=a_j,
 \qquad
 \int_{\Sph}h\overline g\,d\sigma=\sum_{j\in J}a_jb_j.
\]
Consequently \eqref{eq:dual-pairing-estimate} implies
\begin{equation}\label{eq:dual-finite-sequence}
 \sum_{j\in J}a_jb_j
 \le C\norm f_{4/3}\Big(\sum_{j\in J}a_j^4\Big)^{1/4}.
\end{equation}
This is precisely the duality between $\ell^4$ and $\ell^{4/3}$.
For completeness, if $\sum_{j\in J}b_j^{4/3}>0$, choose
\[
 a_j=\frac{b_j^{1/3}}
 {\bigl(\sum_{i\in J}b_i^{4/3}\bigr)^{1/4}}.
\]
Then $\sum a_j^4=1$ and the left side of
\eqref{eq:dual-finite-sequence} equals
$\bigl(\sum_{j\in J}b_j^{4/3}\bigr)^{3/4}$.
We conclude that
\[
 \Big(\sum_{j\in J}\norm{D_jh}_2^{4/3}\Big)^{3/4}
 \le C\norm f_{4/3},
\]
with a constant independent of $J$. Letting $J$ exhaust all dyadic
indices and applying monotone convergence gives
\begin{equation}\label{eq:dual-sharp-result}
 \Big(\sum_{j\ge-1}\norm{D_jh}_2^{4/3}\Big)^{3/4}
 \le C\norm f_{4/3}.
\end{equation}
Equation \eqref{eq:dual-sharp-equivalence} with $q=4/3$ is exactly
\eqref{eq:dual} for Schwartz functions.

\smallskip\noindent\emph{Extension to $L^{4/3}$ and identification.}
For arbitrary $f\in L^{4/3}(\RR^3)$ choose
$f_n\in C_c^\infty(\RR^3)$ with $f_n\to f$ in $L^{4/3}$.
Applying \eqref{eq:dual} to $f_n-f_m$ shows that
$\widehat f_n|_{\Sph}$ is Cauchy in the complete space
$B^0_{2,4/3}(\Sph)$. Denote its limit by $Rf$. The estimate
passes to the limit and gives
$\norm{Rf}_{B^0_{2,4/3}}\le C\norm f_{4/3}$.
The limit is independent of the approximating sequence.
Furthermore, for every finite spherical-harmonic sum $g$,
\[
 \pair{Rf}{g}_{\Sph}
 =\lim_{n\to\infty}\int_{\RR^3}f_n(x)\overline{Eg(x)}\,dx
 =\int_{\RR^3}f(x)\overline{Eg(x)}\,dx,
\]
where the last limit follows from $Eg\in L^4$.
This proves that $Rf$ is the Fourier restriction of $f$ in the
natural distributional sense. If $f\in L^1\cap L^{4/3}$,
choose the approximants to converge in both $L^1$ and $L^{4/3}$.
Their Fourier transforms then converge uniformly, and hence the Besov
limit $Rf$ agrees with the pointwise restriction of the continuous
Fourier transform of $f$.
The identification and the desired bound are proved.
\end{proof}

\subsection{Organization of the proof}
The proof of Theorem~\ref{thm:bilinear} separates physical space according to the radial variable $r=|x|$. In the inner region $r<K/10$, spherical Bessel functions of degree approximately $K$ are exponentially small. At radii $r\ge L^2$, the low band has an exact finite spherical-Hankel expansion; a radial Plancherel identity and a local mass estimate yield the desired $L/K$ factor.

The remaining interval, $K/10\lesssim r<L^2$, is the Fresnel region. Here local angular momentum $r\sin\ang(\omega,\widehat x)$ distinguishes low- and high-frequency rays. The angular separation suggests a cap estimate with factor $L/K$, but the original data are not cap-supported. An exact angular-momentum differential identity supplies above-band and below-band leakage bounds. After localization, a circle-fiber convolution estimate controls the main transverse interaction. A tube-energy bound controls leakage against the transverse high band, while weighted sector counting prevents a logarithmic accumulation over radial shells. Each of these mechanisms is proved separately below, so the full argument does not require treating the phase-space heuristic as an actual support statement.

\section{Spectral reduction and local estimates}\label{sec:spectral}
The plane-wave expansion reads
\begin{equation}\label{eq:pw}
 Eg(r\theta)=4\pi\sum_{k=0}^{\infty}i^k j_k(r)(\Pi_kg)(\theta),
 \quad\theta\in\Sph,\ r>0,
\end{equation}
where $j_k$ is the spherical Bessel function and $\Pi_k$ is the orthogonal projector onto $\HH_k$. In particular, $P_j$ commutes with extension at each fixed radius:
\begin{equation}\label{eq:intertwining}
 P_j(Eg(r\cdot))=E(P_jg)(r\cdot).
\end{equation}
Spherical Littlewood--Paley theory, applied at every radius and integrated with $r^2\,dr$, yields
\begin{equation}\label{eq:LP4}
 \norm{Eg}_4^4\le C\int_{\RR^3}\left(\sum_j|E(P_jg)(x)|^2\right)^2dx
 =C\sum_{j,k}\norm{E(P_jg)E(P_kg)}_2^2.
\end{equation}
We spell out why the double sum is controlled by the \emph{fourth}
powers of the input block norms. Put $b_j=\norm{P_jg}_2^2$.
There is a fixed integer $d_0$ depending only on the separation
threshold $C_0$ and the finite width of the smooth blocks such that,
when $|j-k|\le d_0$, H\"older and Stein--Tomas give
\[
 \norm{E(P_jg)E(P_kg)}_2^2\le Cb_jb_k.
\]
When $|j-k|>d_0$, the separated-band theorem gives
$C2^{-|j-k|}b_jb_k$, after absorbing a fixed factor from the
spectral cutoffs. Define the even nonnegative sequence
\[
 w_d=\begin{cases}
 1,& |d|\le d_0,\\
 2^{-|d|},& |d|>d_0.
 \end{cases}
\]
It has finite $\ell^1(\mathbb Z)$ norm. Thus \eqref{eq:LP4}
becomes
\[
 \norm{Eg}_4^4\le C\sum_{j,k}w_{j-k}b_jb_k
      =C\sum_jb_j(w*b)_j.
\]
Cauchy--Schwarz in $j$, followed by Young's convolution inequality,
gives
\begin{equation}\label{eq:discrete}
 \sum_{j,k}w_{j-k}b_jb_k
 \le\norm b_{\ell^2}\norm{w*b}_{\ell^2}
 \le\norm w_{\ell^1}\norm b_{\ell^2}^2
 =\norm w_{\ell^1}\sum_j\norm{P_jg}_2^4.
\end{equation}
This proves Theorem~\ref{thm:main} for finite angular sums after
Theorem~\ref{thm:bilinear} has been established. Completing such
sums in $B^0_{2,4}$ extends the bound to all data in the stated
space, exactly as asserted in Theorem~\ref{thm:main}.
More generally, any uniform off-diagonal factor
$(L/K)^\varepsilon$ with fixed $\varepsilon>0$ would give a
summable sequence $w$ and suffice for the Besov conclusion; the
present paper proves the sharper $L/K$ decay.

We fix the following standard estimates, stating explicitly the forms used below.
\begin{lemma}\label{lem:kernel}
For a sufficiently large fixed integer $n_0$ and a smooth multiplier $a$ supported in $[1/8,4]$, its kernel $K_{a(\Lambda/L)}$ satisfies, for $L\ge1$,
\begin{equation}\label{eq:needlet}
 |K_{a(\Lambda/L)}(\omega,\eta)|\le C\norm{a}_{C^{n_0}}L^2(1+Ld(\omega,\eta))^{-7}.
\end{equation}
Here $d$ denotes geodesic distance on $\Sph$. The constants are independent of $L$ and $a$.
\end{lemma}
\begin{proof}
The localization theorem for smooth spherical spectral multipliers
\cite{npwSIAM2006,npwJFA2006} gives, for every fixed
$a\in C_c^\infty((1/8,4))$, a constant $C_a$ such that
\[
 \sup_{L\ge1}\sup_{\omega,\eta\in\Sph}
 L^{-2}(1+Ld(\omega,\eta))^7
 |K_{a(\Lambda/L)}(\omega,\eta)|\le C_a.
\]
We explain why the constant can be made uniform in $a$ with a finite
number of derivatives, since this quantitative dependence is used when
$a$ varies with the power of $Q$ later in the proof.
Work in the Fr\'echet space $\mathcal D_{[1/8,4]}$ of smooth functions
supported in the fixed compact interval $[1/8,4]$, endowed with the
seminorms $\|a\|_{C^j}$. For each choice of $L,\omega,\eta$, the
linear functional
\[
 a\longmapsto L^{-2}(1+Ld(\omega,\eta))^7
 K_{a(\Lambda/L)}(\omega,\eta)
\]
 is continuous: at each $L$ the spectral sum involves only finitely
many eigenvalues. The displayed localization bound says that this
family of functionals is pointwise bounded on
$\mathcal D_{[1/8,4]}$. By the uniform boundedness principle for
Fr\'echet spaces, the family is equicontinuous. Since its topology
is defined by the increasing seminorms $\|a\|_{C^j}$, there are a
fixed integer $n_0$ and a constant $C$, independent of $L$, $\omega$,
$\eta$, and $a$, such that the family is bounded by
$C\|a\|_{C^{n_0}}$. This proves \eqref{eq:needlet}.
\end{proof}
\begin{lemma}\label{lem:localmass}
For $v\in\HH_k$ and every geodesic cap $B(\zeta,\rho)$, $0<\rho\le1$,
\begin{equation}\label{eq:capmass}
 \int_{B(\zeta,\rho)}|v|^2\,d\sigma\le C\rho\norm{v}_2^2.
\end{equation}
\end{lemma}
\begin{proof}
If $\rho\ge(k+1)^{-1}$, this is the uniform local $L^2$-mass estimate on the sphere \cite{sogge2016}. For $\rho<(k+1)^{-1}$, the addition theorem gives
$\norm{v}_\infty^2\le(2k+1)(4\pi)^{-1}\norm{v}_2^2$; multiplying by $\sigma(B)\le\pi\rho^2$ gives the result. The case $k=0$ is immediate.
\end{proof}
\begin{lemma}\label{lem:static}
If $u$ has $\Lambda$-spectrum in $[L/2,2L]$ and $v\in\HH_k$ is arbitrary, then
\begin{equation}\label{eq:static}
 \int_{\Sph}|u|^2|v|^2\,d\sigma\le CL\norm{u}_2^2\norm{v}_2^2.
\end{equation}
\end{lemma}
\begin{proof}
Choose $\phi\in C_c^\infty((1/8,4))$ equal to one on $[1/4,3]$, so $u=\phi(\Lambda/L)u$. By Lemma \ref{lem:kernel}, Cauchy--Schwarz, and the fact that $L^2(1+Ld)^{-7}$ has uniformly bounded integral,
\begin{equation}\label{eq:staticpointwise}
 |u(\theta)|^2\le CL^2\int_{\Sph}(1+Ld(\theta,\eta))^{-7}|u(\eta)|^2\,d\sigma(\eta).
\end{equation}
We now explain the step from the pointwise kernel to the $L$ factor; the latter must be uniform in the much higher degree of $v$. Fix $\eta$. In the first cap $d(\theta,\eta)<L^{-1}$, Lemma~\ref{lem:localmass} bounds the $|v|^2$ mass by $CL^{-1}\norm v_2^2$. In the annulus with $2^{a-1}L^{-1}\le d(\theta,\eta)<2^aL^{-1}$, the kernel loses a factor at most $C2^{-7a}$ and the same lemma bounds the mass by $C(2^a/L)\norm v_2^2$. Consequently the sum is dominated by
$CL^{-1}\sum_{a\ge0}2^{-6a}\norm v_2^2$. This gives, with all radial scales made explicit,
\begin{equation*}
 \sup_\eta\int_{\Sph}(1+Ld(\theta,\eta))^{-7}|v(\theta)|^2\,d\sigma(\theta)
 \le CL^{-1}\norm{v}_2^2.
\end{equation*}
Insert this into \eqref{eq:staticpointwise} and integrate in $\theta$.
\end{proof}

\section{Circle fibers, radial structure, and tubes}
\begin{lemma}\label{lem:fiber}
Let $f_1$ be supported in the two caps $B(e,\delta)\cup B(-e,\delta)$, and suppose
$|\omega\times\eta|\ge\vartheta>0$ for all $\omega\in\supp f_1$, $\eta\in\supp f_2$. Then
\begin{equation}\label{eq:fiber}
 \norm{Ef_1Ef_2}_{L^2(\RR^3)}^2\le C\frac\delta\vartheta\norm{f_1}_2^2\norm{f_2}_2^2.
\end{equation}
\end{lemma}
\begin{proof}
Plancherel gives $\norm{Ef_1Ef_2}_2^2=(2\pi)^3\norm{f_1d\sigma*f_2d\sigma}_2^2$. For $0<|\xi|<2$, the fiber $\{\omega\in\Sph:|\xi-\omega|=1\}$ is a circle $C_\xi$ in a plane orthogonal to $\xi$. Coarea gives
\begin{equation}\label{eq:coarea}
 (f_1d\sigma*f_2d\sigma)(\xi)=\int_{C_\xi}
 \frac{f_1(\omega)f_2(\xi-\omega)}{|\omega\times(\xi-\omega)|}\,d\HH^1(\omega).
\end{equation}
Let $G_\xi$ be the portion of $C_\xi$ on which both factors in \eqref{eq:coarea} are nonzero, and write $d\mu_\xi=d\mathcal H^1/|\omega\times(\xi-\omega)|$. For every $\xi$, Cauchy--Schwarz yields
\[
 |(f_1d\sigma*f_2d\sigma)(\xi)|^2
 \le\mu_\xi(G_\xi)
 \int_{G_\xi}|f_1(\omega)|^2|f_2(\xi-\omega)|^2\,d\mu_\xi(\omega).
\]
The cross-product hypothesis bounds the Jacobian denominator below by $\vartheta$. Since the first datum occupies two caps of radius $\delta$, each circle intersects these caps in total arc length at most $C\delta$, whether the circle has radius smaller or larger than $\delta$. Thus $\mu_\xi(G_\xi)\le C\delta/\vartheta$. Integrating in $\xi$ and applying the same coarea formula to the remaining integral gives $\norm{f_1}_2^2\norm{f_2}_2^2$, proving \eqref{eq:fiber}. Separation from both parallel and antiparallel directions is needed because the convolution Jacobian vanishes in either configuration.
\end{proof}

\begin{lemma}\label{lem:ray}
For $g_K$ supported in $\Lambda\ge K/2$,
\begin{equation}\label{eq:ray}
 \int_0^\infty|Eg_K(r\theta)|^2\,dr
 \le\frac{16\pi^3}{K}\sum_k|(\Pi_kg_K)(\theta)|^2.
\end{equation}
Furthermore, with
\begin{equation}\label{eq:a_n}
 a_n(k)=\begin{cases}\displaystyle\frac{(k+n)!}{n!(k-n)!2^n},&0\le n\le k,\\0,&n>k,\end{cases}
\end{equation}
let $u_n=\sum_ka_n(k)\Pi_kg$. Then the exact identity
\begin{equation}\label{eq:hankel}
 rEg(r\theta)=2\pi\sum_{n\ge0}r^{-n}
 \left(-i^{n+1}e^{ir}u_n(\theta)+i(-i)^ne^{-ir}u_n(-\theta)\right)
\end{equation}
holds for finite harmonic sums. For $g_L$ with degrees $k\le2L$,
\begin{equation}\label{eq:un}
 \norm{u_n}_2\le \frac{(2L^2)^n}{n!}\norm{g_L}_2.
\end{equation}
\end{lemma}
\begin{proof}
For completeness, define $H_\theta(s)=\sum_kP_k(s)(\Pi_kg)(\theta)$ for $-1\le s\le1$, extended by zero outside this interval. Then the plane-wave identity gives $Eg(r\theta)=2\pi\int_{-1}^1e^{irs}H_\theta(s)\,ds$. Plancherel in the one-dimensional variable $r$ and orthogonality $\int_{-1}^1P_k(s)P_{k'}(s)\,ds=2\delta_{kk'}/(2k+1)$ yield
\begin{equation*}
 \int_{\RR}|Eg(r\theta)|^2\,dr=16\pi^3\sum_k\frac{|(\Pi_kg)(\theta)|^2}{2k+1},
\end{equation*}
from orthogonality of the Legendre polynomials. Restricting to $r\ge0$ proves \eqref{eq:ray}. The spherical Hankel functions have a finite polynomial expansion in $r^{-1}$ at each integral order $k$, with coefficients $a_n(k)$ displayed in \eqref{eq:a_n}; see \cite[\S10.49]{dlmf}. Writing $j_k=(h_k^{(1)}+h_k^{(2)})/2$, inserting the expansions into \eqref{eq:pw}, and using $Y_k(-\theta)=(-1)^kY_k(\theta)$ leads to the two outgoing/incoming terms in \eqref{eq:hankel}. Because the sum in $n$ is finite for each harmonic, no interchange of infinite asymptotic series is involved for finite input. Finally, $a_n(k)\le [k(k+1)/2]^n/n!$ and $k+1/2\le2L$ imply \eqref{eq:un}.
\end{proof}
\begin{lemma}\label{lem:tube}
Let $e\in\Sph$ and $f\in L^2(\Sph)$ have support in $\{\eta:\sin\ang(\eta,e)\ge\beta\}$, $\beta>0$. For $W>0$,
\begin{equation}\label{eq:tube}
 \int_{\dist(x,\RR e)\le W}|Ef(x)|^2\,dx\le C\frac W\beta\norm{f}_2^2.
\end{equation}
Consequently, for $R\ge1$,
\begin{equation}\label{eq:AH}
 \int_{|x|\le R}|Ef(x)|^2\,dx\le CR\norm{f}_2^2.
\end{equation}
\end{lemma}
\begin{proof}
Write $\eta=ue+(1-u^2)^{1/2}\vartheta$ with $\vartheta\in S^1\subset e^\perp$, and $x=se+y$, $y\in e^\perp$. Archimedes' formula gives $d\sigma=du\,d\vartheta$. Plancherel in $s$ reduces the tube integral to
\begin{equation*}
 2\pi\int_{-1}^1du\int_{|y|\le W}
 \left|\int_{S^1}e^{i(1-u^2)^{1/2}y\cdot\vartheta}f(u,\vartheta)\,d\vartheta\right|^2dy.
\end{equation*}
For an unscaled circle, Agmon--H\"ormander's local estimate gives
$\int_{|z|\le T}|\widehat{h\,d\vartheta}(z)|^2\,dz\le CT\norm h_2^2$ when $T\ge1$. For $T<1$, Cauchy--Schwarz gives the even stronger $CT^2\norm h_2^2$ bound. Substitute $z=(1-u^2)^{1/2}y$ in the preceding display: the Jacobian contributes $(1-u^2)^{-1}$, while the circle bound contributes $(1-u^2)^{1/2}W$. Their product is $CW/(1-u^2)^{1/2}$, at most $CW/\beta$. Integrating in $u$ establishes the tube inequality.

For the ball inequality \eqref{eq:AH}, choose two orthogonal axes $e_1,e_2$. Every direction on $\Sph$ has sine of its angle with at least one of these axes bounded below by a fixed positive constant. A smooth two-piece partition therefore splits the datum into functions to which the tube inequality applies. Since $B_R$ lies in the width-$R$ tube around either axis, summing the two estimates gives \eqref{eq:AH}.
\end{proof}

\section{Angular-momentum leakage away from the characteristic angles}
This is the main technical part of the argument. A heuristic phase-space picture says that harmonic degree $N$ is associated with angular momentum $|x\times\omega|\simeq N$. A heuristic cannot by itself justify cutting a spherical-harmonic band into angular caps: such a cutoff changes its spectrum. We therefore prove a leakage estimate directly, without assuming that the original data had cap support.

The key estimate is a quantitative form of the fact that the plane wave $\exp(-ix\cdot\omega)$ has local angular momentum of magnitude $|x\times\omega|$. The gain is expressed in terms of the Fresnel number $|x|\psi^2$; a bound involving only the width/separation ratio would not suffice uniformly in our final radial sum.

Write
\begin{equation}\label{eq:angmomentum}
 \LL=-i\omega\times\nabla_\omega,
 \qquad Q=\LL\cdot\LL,
 \qquad \mathfrak e_x(\omega)=e^{-ix\cdot\omega}.
\end{equation}
For $x\ne0$, let $t=|x|$, $v=x\times\omega$, $p=x\cdot\omega$, and $s=|v|^2$. If $\psi\le\pi/4$ and $f$ is supported in an angular region at scale $\psi$, introduce
\begin{equation}\label{eq:seminorm}
 [f]_{n,\psi}=\max_{|\alpha|\le n}\psi^{|\alpha|}\norm{\LL^\alpha f}_{L^\infty(\Sph)},
\end{equation}
where $\alpha$ is a word in the three components of $\LL$. All angular cutoffs below have uniformly bounded $[\cdot]_{n,\psi}$ at their natural scale.

\begin{lemma}\label{lem:leak}
Fix an integer $m\ge1$. There exist $A_m\ge1$, $0<a_m<1$ and $C_m<\infty$ with the following properties. Let $t=|x|$, and $0<\psi\le\pi/4$.
\begin{enumerate}[label=\textup{(\alph*)}]
\item \emph{Above-band leakage.} Suppose $u$ has $Q$-spectrum in $[0,N^2]$, and $\chi$ is supported where
\begin{equation}\label{eq:upperregion}
 \ang(\omega,\RR\widehat x)\ge\psi,
 \qquad |x\times\omega|\le64t\psi.
\end{equation}
If $t\psi\ge A_mN$ and $t\psi^2\ge1$, then
\begin{equation}\label{eq:above}
 |E(\chi u)(x)|\le C_m(t\psi^2)^{-m}[\chi]_{2m,\psi}
 |\supp\chi|^{1/2}\sum_{i\ge0}2^{-i}N^{-2i}
 \norm{Q^iu}_{L^2(\supp\chi)}.
\end{equation}
\item \emph{Below-band leakage.} Suppose $u$ has $Q$-spectrum in $[N^2,\infty)$, and $\chi$ is supported in the union of the caps of radius $\psi$ around $\widehat x$ and $-\widehat x$. If $t\psi\le a_mN$, then
\begin{equation}\label{eq:below}
 |E(\chi u)(x)|\le C_m\left(\frac{t+\psi^{-2}}{N^2}\right)^m
 [\chi]_{2m,\psi}|\supp\chi|^{1/2}\norm{u}_2.
\end{equation}
\end{enumerate}
The estimates extend to endpoint regions reaching the equator, provided the upper bound $|x\times\omega|\le64t\psi$ is retained.
\end{lemma}
\begin{proof}
Direct differentiation gives
\begin{equation}\label{eq:Le}
 \LL_a\mathfrak e_x=v_a\mathfrak e_x,
 \quad \sum_{a=1}^3\LL_av_a=-2ip,
\end{equation}
and hence the exact angular-momentum identity
\begin{equation}\label{eq:Qe}
 Q(f\mathfrak e_x)=(sf+\mathcal Rf)\mathfrak e_x,
 \qquad \mathcal Rf=Qf+2v\cdot\LL f-2ipf.
\end{equation}
The components $v_a$ and $p$ are linear functions of $\omega$ with coefficients bounded by $t$. Differentiating such functions under the seminorm \eqref{eq:seminorm} shows, on a region where $|v|\le64t\psi$,
\begin{equation}\label{eq:symbolbounds}
 [sf]_{n,\psi}\le C_nt^2\psi^2[f]_{n,\psi},\qquad
 [\mathcal Rf]_{n,\psi}\le C_n(t+\psi^{-2})[f]_{n+2,\psi}.
\end{equation}
If in addition the angular distance from $\pm\widehat x$ is at least $\psi$, then $s\gtrsim t^2\psi^2$ and the same weighted product-rule calculation gives
\begin{equation}\label{eq:sinverse}
 [f/s]_{n,\psi}\le C_n(t\psi)^{-2}[f]_{n,\psi}.
\end{equation}
In each case the differential operators preserve support. Here and below the constants can be increased once for the finitely many derivative orders at most $2m$.

For (a), self-adjointness of $Q$ and \eqref{eq:Qe} give, for $T(w,f)=\pair{w}{f\mathfrak e_x}$,
\begin{equation}\label{eq:aboverec}
 T(w,f)=T(Qw,f/s)-T(w,\mathcal R(f/s)).
\end{equation}
Call the first term a $B$-step and the second a $G$-step. By \eqref{eq:symbolbounds}--\eqref{eq:sinverse}, a $B$-step replaces $w$ by $Qw$ and costs at most $C(t\psi)^{-2}$ without increasing the derivative order required of $f$. A $G$-step costs at most
\begin{equation}\label{eq:Gcost}
 C\frac{t+\psi^{-2}}{t^2\psi^2}\le\frac{C'}{t\psi^2},
\end{equation}
and increases the derivative order required of $f$ by two; the last inequality uses $t\psi^2\ge1$. Iterate \eqref{eq:aboverec}, stopping a branch after its $m$th $G$-step. Each completed word with $i$ $B$-steps has at most $m$ $G$-steps and is bounded by
\begin{equation}\label{eq:word}
 C_m\bigl(C_m(t\psi)^{-2}\bigr)^i(t\psi^2)^{-m}
 [\chi]_{2m,\psi}|\supp\chi|^{1/2}
 \norm{Q^iu}_{L^2(\supp\chi)}.
\end{equation}
The number of words with exactly $i$ principal steps before the $m$-th remainder step is $\binom{i+m-1}{m-1}$. This grows polynomially in $i$ when $m$ is fixed, so the convergence condition is genuinely a geometric one and cannot be sacrificed to the number of branches. Choose $A_m$ so large that the coefficient $C_mN^2/(t\psi)^2$ is smaller than a sufficiently small constant, say $1/4$ after absorbing the finitely many derivative constants. Then the polynomial word count times this coefficient to the power $i$ is bounded by $C_m2^{-i}$, and the sum of completed words is bounded by the right-hand side of \eqref{eq:above}.

For rigor one must also account for words that have not yet accumulated $m$ remainder steps after $n$ iterations. Such a word contains at least $n-m+1$ principal steps and at most $m-1$ remainder steps. There are $O_m(n^{m-1})$ possible words. The spectral inequality $\norm{Q^ju}_2\le N^{2j}\norm u_2$ and the same geometric smallness show that their total tends to zero as $n\to\infty$. We may therefore take the limit in the finite expansion, which establishes \eqref{eq:above} without leaving an unestimated remainder.

For (b), the spectral inverse is defined on the range of $u$ and satisfies $\norm{Q^{-1}w}_2\le N^{-2}\norm w_2$. Equation \eqref{eq:Qe} yields instead
\begin{equation}\label{eq:belowrec}
 T(w,f)=T(Q^{-1}w,sf)+T(Q^{-1}w,\mathcal Rf).
\end{equation}
The first step costs $Ct^2\psi^2/N^2$ and the second costs $C(t+\psi^{-2})/N^2$ and two derivatives of the cutoff. This time the angular momentum of the phase is too small for the spectral support of $u$. Choose $a_m$ so that the first cost is at most $1/4$. Expand the identity to depth $n$, collect words after the $m$-th remainder step, and note that each word has at most $2m$ derivatives of the cutoff. As above, the polynomial number of paths is dominated by the geometric principal-step factor. The unfinished paths tend to zero by the bound $\norm{Q^{-j}u}_2\le N^{-2j}\norm u_2$. Taking the limit leaves precisely $m$ powers of $(t+\psi^{-2})/N^2$, which proves \eqref{eq:below}.
\end{proof}

The local norm in \eqref{eq:above} must ultimately be controlled by one input norm without paying for the total number of sectors. The following lemma supplies that control.
\begin{lemma}\label{lem:localpowers}
For $L\ge1$ and measurable $U\subset\Sph$, set
\begin{equation}\label{eq:weightnorm}
 \cN_U(g)^2=\int_{\Sph}|g(\omega)|^2(1+Ld(\omega,U))^{-4}\,d\sigma(\omega).
\end{equation}
If $g_L$ has $\Lambda$-spectrum in $[L/2,2L]$, then, for every integer $i\ge0$,
\begin{equation}\label{eq:localpowers}
 \norm{Q^ig_L}_{L^2(U)}\le C\,17^iL^{2i}\cN_U(g_L).
\end{equation}
\end{lemma}
\begin{proof}
Take $\phi\in C_c^\infty((1/8,4))$ equal to one on $[1/2,2]$. Since $g_L=\phi(\Lambda/L)g_L$,
\begin{equation*}
 Q^ig_L=L^{2i}\phi_i(\Lambda/L)g_L,
 \qquad \phi_i(z)=\bigl(z^2-(2L)^{-2}\bigr)^i\phi(z).
\end{equation*}
For the fixed derivative order $n_0$ of Lemma \ref{lem:kernel}, $\norm{\phi_i}_{C^{n_0}}\le C(1+i)^{n_0}16^i\le C'17^i$, with $C'$ independent of $L,i$. The kernel estimate \eqref{eq:needlet}, together with Cauchy--Schwarz and the triangle inequality for geodesic distance, implies the weighted Schur estimate
\begin{equation*}
 \int_U\left|\int L^2(1+Ld(\omega,\eta))^{-7}g(\eta)\,d\sigma(\eta)\right|^2d\sigma(\omega)
 \le C\int_{\Sph}|g(\eta)|^2(1+Ld(\eta,U))^{-4}d\sigma(\eta).
\end{equation*}
Indeed, split the kernel into the factor $(1+Ld)^{-3}$, whose $L^2$-scaled integral is bounded, and the remaining off-diagonal factor; use $d(\eta,U)\le d(\eta,\omega)$ for $\omega\in U$. Substituting the bound for $\phi_i$ proves \eqref{eq:localpowers}.
\end{proof}
\begin{remark}\label{rem:threshold}
The use of $N=4L$ in Lemma \ref{lem:leak}(a) gives
\begin{equation}\label{eq:17over32}
 \sum_{i\ge0}2^{-i}(4L)^{-2i}\norm{Q^ig_L}_{L^2(U)}
 \le C\sum_{i\ge0}(17/32)^i\cN_U(g_L)\le C\cN_U(g_L).
\end{equation}
The explicit geometric ratio is why the local derivative norm does not introduce a hidden power of $L$.
\end{remark}

\section{Inner and far-field contributions}
We prove Theorem \ref{thm:bilinear}. By density it is enough to consider finite harmonic sums, with spectra as in \eqref{eq:bands}. Fix a large separation constant $C_0$ later and suppose $K\ge C_0L$.

\subsection{The inner ball}
On $|x|\le K/10$, the degrees of $g_K$ satisfy $k\ge K/2-1/2$. The elementary bound
\begin{equation*}
 |j_k(r)|\le\frac{r^k}{(2k+1)!!}\le \left(\frac{er}{2k}\right)^k
\end{equation*}
is exponentially small when $r\le K/10$. The addition theorem and \eqref{eq:pw} give $|Eg_K(x)|\le Ce^{-cK}\norm{g_K}_2$, while Cauchy--Schwarz directly in \eqref{eq:Edef} gives $|Eg_L(x)|\le(4\pi)^{1/2}\norm{g_L}_2$. Thus
\begin{equation}\label{eq:inner}
 \int_{|x|<K/10}|Eg_L|^2|Eg_K|^2\,dx
 \le CK^3e^{-2cK}\norm{g_L}_2^2\norm{g_K}_2^2
 \le C\frac LK\norm{g_L}_2^2\norm{g_K}_2^2.
\end{equation}

\subsection{The far field}
Let
\begin{equation}\label{eq:R0}
 R_0=\max\{L^2,K/10\}.
\end{equation}
From \eqref{eq:hankel}, Cauchy--Schwarz with the weights $2^{-n}$ gives
\begin{equation}\label{eq:farpoint}
 r^2|Eg_L(r\theta)|^2
 \le C\sum_{n\ge0}2^nr^{-2n}\bigl(|u_n(\theta)|^2+|u_n(-\theta)|^2\bigr),\quad r\ge R_0.
\end{equation}
Integrate $r\ge R_0$ against $|Eg_K(r\theta)|^2$ and then in $\theta$. Apply Lemma \ref{lem:ray}, the static product estimate \eqref{eq:static} to each $u_n$ and every $\Pi_kg_K$, and \eqref{eq:un}. We obtain
\begin{align}
 \int_{|x|\ge R_0}|Eg_L|^2|Eg_K|^2\,dx
 &\le C\frac LK\norm{g_K}_2^2
 \sum_{n\ge0}2^nR_0^{-2n}\norm{u_n}_2^2\notag\\
 &\le C\frac LK\left(\sum_{n\ge0}\frac{8^n}{(n!)^2}\right)
 \norm{g_L}_2^2\norm{g_K}_2^2.\label{eq:far}
\end{align}
This proves the required estimate for the far field, independently of the Fresnel localization argument.

\section{The Fresnel range}
Only the region
\begin{equation}\label{eq:Fresnel}
 K/10\le |x|<L^2
\end{equation}
remains. It is empty when $K\ge10L^2$. We now fix all constants in an order that avoids circular dependencies. Choose $m=2$ in Lemma \ref{lem:leak}, and let $A_2,a_2,C_2$ be its corresponding constants. Choose
\begin{equation}\label{eq:constants}
 \kappa\ge\max\{60,8A_2,6\sqrt5/a_2\},\quad
 M\ge\max\{2,16A_2\},\quad C_0\ge4\kappa M,
\end{equation}
with each constant now fixed absolutely. All implicit constants in this section may depend on $\kappa$ and $M$, but not on $L,K$ or the input functions.

\subsection{Sectors and cutoffs}
Decompose \eqref{eq:Fresnel} into dyadic shells $S_R=\{R\le|x|<2R\}$ with $R\simeq2^jK/10<L^2$. For a given shell put $\delta=L/R$. Thus
\begin{equation}\label{eq:scales}
 L\delta=L^2/R>1,
 \quad \delta\le10L/K\ll1,
 \quad R\ge K/10.
\end{equation}
Choose a maximal $\delta$-separated set of directions $\{\omega_c\}$ on the sphere, and define
\begin{equation}\label{eq:sector}
 \Gamma_{c,R}=\{x\in S_R:\ang(\widehat x,\omega_c)\le\delta\}.
\end{equation}
They cover $S_R$ with bounded overlap. The low-frequency cutoff $\chi_c$ equals one for $\ang(\omega,\RR\omega_c)\le M\delta$ and is supported where that angle is at most $2M\delta$. The high-frequency cutoff $\widetilde\chi_c$ equals one on the angular annulus
\begin{equation}\label{eq:highcut}
 \frac{2K}{\kappa R}\le\ang(\omega,\RR\omega_c)\le\frac{\kappa K}{2R}
\end{equation}
and is supported in the enlarged annulus with lower boundary $K/(\kappa R)$ and upper boundary $\kappa K/R$. Angles are restricted, of course, to $[0,\pi/2]$, and the high cutoff can equal one all the way to the equator if its upper boundary exceeds $\pi/2$. Both cutoffs are smooth, even under $\omega\mapsto-\omega$, and have bounded scale-normalized derivatives of every fixed order.

Set
\begin{equation}\label{eq:splits}
 f_c=\chi_cg_L,\quad h_c=\widetilde\chi_cg_K,
 \quad \ell_c=E((1-\chi_c)g_L),\quad q_c=E((1-\widetilde\chi_c)g_K).
\end{equation}
For $x\in\Gamma_{c,R}$, expand the squared product and use
\begin{equation}\label{eq:four}
 |Eg_L(x)Eg_K(x)|^2\le4\bigl(I_{c,R}+II_{c,R}+III_{c,R}+IV_{c,R}\bigr),
\end{equation}
where pointwise $I=|Ef_cEh_c|^2$, $II=|\ell_cEh_c|^2$, $III=|Ef_cq_c|^2$, $IV=|\ell_cq_c|^2$. We estimate the integrals of these terms over the sectors and sum.

\subsection{The transverse main term}
If $\omega\in\supp f_c$ and $\eta\in\supp h_c$, the inequality $K\ge4\kappa ML$ gives angular separation
\begin{equation}\label{eq:transverse}
 |\omega\times\eta|\ge c\frac K{\kappa R}.
\end{equation}
Lemma \ref{lem:fiber} with cap radius $2M\delta$ therefore gives
\begin{equation}\label{eq:mainsector}
 \int_{\Gamma_{c,R}}I_{c,R}\,dx\le C\kappa M\frac LK\norm{f_c}_2^2\norm{h_c}_2^2.
\end{equation}
For a fixed pair $(\omega,\eta)$, only $O(M^2)$ sector centers $\omega_c$ can lie within $2M\delta$ of $\pm\omega$. Moreover, the angle $\alpha=\ang(\omega,\RR\eta)$ then belongs to an interval comparable to $[K/(\kappa R),\kappa K/R]$. Hence a fixed pair contributes for at most $O(\log\kappa)$ dyadic radii, uniformly in $K,L$. Fubini and the bounded overlap count give
\begin{equation}\label{eq:mainall}
 \sum_{R,c}\int_{\Gamma_{c,R}}I_{c,R}\,dx\le C\kappa M^3\log(2\kappa)\frac LK
 \norm{g_L}_2^2\norm{g_K}_2^2.
\end{equation}
The absence of a factor $\log(L^2/K)$ here is essential. It comes from a fact that is easy to overlook: fixing a pair of input directions also fixes the characteristic radial scale up to a constant number of dyadic choices. We are not simply adding the same cap estimate over every possible physical radius.

\subsection{Low-band leakage against the transverse high band}
For every $c,R$, decompose $1-\chi_c$ into a smooth dyadic partition $\sum_j\chi_{c,j}$ on angular annuli outside the low cap. For the ordinary annuli choose
\begin{equation}\label{eq:annuli}
 \supp\chi_{c,j}\subset\bigl\{\ang(\omega,\RR\omega_c)\in[2^jM\delta,2^{j+2}M\delta]\bigr\},
 \qquad \psi_j=2^{j-1}M\delta.
\end{equation}
Stop when the annuli approach the equator; the final piece is supported on the remaining angular region. Choosing the last index $J$ by $2^{J+5}M\delta\le\pi/2<2^{J+6}M\delta$ ensures that its associated scale $\psi_{J+1}=2^JM\delta$ lies between $\pi/128$ and $\pi/64$. Consequently all pieces, including the last, satisfy the support condition \eqref{eq:upperregion} relative to any $x\in\Gamma_{c,R}$, after a harmless fixed enlargement of the angular annuli. For each $j$ let $U_{c,j}$ be such a fixed enlarged annulus around $\pm\omega_c$, containing the support of $\chi_{c,j}$. Define
\begin{equation}\label{eq:Ncj}
 \cN_{c,j}=\cN_{U_{c,j}}(g_L)
 =\left(\int_{\Sph}|g_L(\omega)|^2(1+Ld(\omega,U_{c,j}))^{-4}\,d\sigma(\omega)\right)^{1/2}.
\end{equation}

For $x\in\Gamma_{c,R}$, the cutoff $\chi_{c,j}$ has $[\chi_{c,j}]_{4,\psi_j}\le C$, and
\begin{equation}\label{eq:aboveshell}
 t\psi_j\ge\frac M2L\ge A_2(4L),\qquad
 t\psi_j^2\gtrsim4^j\frac{M^2L^2}{R}\ge1.
\end{equation}
Apply Lemma \ref{lem:leak}(a) with $N=4L$ and Lemma \ref{lem:localpowers}, using \eqref{eq:17over32}. Since $|\supp\chi_{c,j}|^{1/2}\le C2^jM\delta$, we get
\begin{equation}\label{eq:lowleakpiece}
 \sup_{x\in\Gamma_{c,R}}|E(\chi_{c,j}g_L)(x)|
 \le C2^{-3j}M\delta\left(\frac R{M^2L^2}\right)^2\cN_{c,j}.
\end{equation}
Here and below, the last annulus is interpreted with the fixed endpoint scale just described. Summation in $j$ and Cauchy--Schwarz give
\begin{equation}\label{eq:lowleak}
 \sup_{x\in\Gamma_{c,R}}|\ell_c(x)|^2
 \le CM^2\delta^2\left(\frac R{M^2L^2}\right)^4
 \sum_j2^{-3j}\cN_{c,j}^2.
\end{equation}

The support of $h_c$ is transverse to the line $\RR\omega_c$, with
$\sin\ang(\omega,\RR\omega_c)\gtrsim K/(\kappa R)$. Every sector $\Gamma_{c,R}$ is contained in the tube $\dist(x,\RR\omega_c)\le2L$. By Lemma \ref{lem:tube},
\begin{equation}\label{eq:transversemass}
 \int_{\Gamma_{c,R}}|Eh_c|^2dx\le C\kappa\frac{LR}K\norm{h_c}_2^2.
\end{equation}
The remaining issue is to sum the weighted local norms over $c$. The $\delta$-separation of the centers and $L\delta>1$ imply
\begin{equation}\label{eq:weightcount}
 \sup_{\omega\in\Sph}\sum_c(1+Ld(\omega,U_{c,j}))^{-4}\le C4^jM^2.
\end{equation}
Indeed, there are $O(4^jM^2)$ centers within distance $O(2^jM\delta)$ of a fixed $\omega$. Beyond that distance, dyadic shells of centers at distance $\rho$ contribute
\begin{equation*}
 C(\rho/\delta)^2(L\rho)^{-4}
 \le C(2^jM)^{-2}(L\delta)^{-4}
\end{equation*}
at the first shell, followed by a convergent geometric series. Since $L\delta>1$, this is bounded. For the terminal annulus one uses the total number of centers, $O(\delta^{-2})$, which is $O(4^{J+1}M^2)$ by the definition of $J$.

Multiply \eqref{eq:lowleak} and \eqref{eq:transversemass}, sum over $c$, and use $\norm{h_c}_2\le\norm{g_K}_2$ together with \eqref{eq:weightcount}. The geometric sum in $j$ gives
\begin{align}
 \sum_c\int_{\Gamma_{c,R}}II_{c,R}\,dx
 &\le C\kappa M^{-4}\frac LK
 \left(\frac R{L^2}\right)^3\norm{g_L}_2^2\norm{g_K}_2^2.
 \label{eq:IIone}
\end{align}
Summing over $R<L^2$ is harmless because $\sum_{R<L^2}(R/L^2)^3\le C$ for dyadic $R$. Thus
\begin{equation}\label{eq:IIall}
 \sum_{R,c}\int_{\Gamma_{c,R}}II_{c,R}\,dx
 \le C\kappa M^{-4}\frac LK\norm{g_L}_2^2\norm{g_K}_2^2.
\end{equation}
The powers in \eqref{eq:IIone} are worth keeping visible. The squared low-band leakage contributes $M^2(L/R)^2(R/(M^2L^2))^4$, the tube energy contributes $LR/K$, and the weighted sector multiplicity contributes $M^2$ after summing the angular annuli. Their product is $M^{-4}(L/K)(R/L^2)^3$ up to the fixed transversality parameter. The positive power of $R/L^2$ makes the remaining dyadic radial sum converge. This is the step that combines the inverse Fresnel number with transverse tube energy to obtain the full factor $L/K$.

\subsection{High-band leakage and the double-error term}
For $x\in\Gamma_{c,R}$, the function $(1-\widetilde\chi_c)g_K$ is supported either near the line $\RR\widehat x$, at angle at most
\begin{equation}\label{eq:psi1}
 \psi_1=3K/(\kappa R),
\end{equation}
or on angular annuli whose distances from that line are bounded below by comparable multiples of $\kappa K/R$. These two portions are separated by smooth cutoffs with uniformly bounded scale-normalized derivatives. We estimate them using the two parts of Lemma \ref{lem:leak}.

For the near part, $Q\ge K^2/5$ on the spectrum of $g_K$ (after enlarging the fixed separation constant if needed). With $N=K/\sqrt5$ we have
\begin{equation}\label{eq:belowconditions}
 t\psi_1\le6K/\kappa\le a_2N,
 \qquad \frac{t+\psi_1^{-2}}{N^2}\le C\frac R{K^2}.
\end{equation}
The second inequality uses $t<2R<L^2+R$, $K\ge C_0L$, and $\psi_1^{-2}=\kappa^2R^2/(9K^2)\le R$. Since the near cutoff has area $O((K/R)^2)$, Lemma \ref{lem:leak}(b) with $m=2$ gives
\begin{equation}\label{eq:nearhigh}
 \sup_{\Gamma_{c,R}}|q_c^{\mathrm{near}}|
 \le C\left(\frac R{K^2}\right)^2\frac KR\norm{g_K}_2
 =C\frac R{K^3}\norm{g_K}_2.
\end{equation}
For the far part, use angular annuli of scale
$\psi_i'\simeq2^i\kappa K/R$, with the final piece allowed to reach the equator. They satisfy $t\psi_i'\ge A_2(2K)$ and $t(\psi_i')^2\ge1$, so Lemma \ref{lem:leak}(a) with $N=2K$ applies. The derivative sum is bounded globally because $Q\le4K^2$ on the band:
\begin{equation*}
 \sum_{j\ge0}2^{-j}(2K)^{-2j}\norm{Q^jg_K}_2\le2\norm{g_K}_2.
\end{equation*}
Using the annulus area and the two inverse-Fresnel factors yields
\begin{equation}\label{eq:farhigh}
 \sup_{\Gamma_{c,R}}|q_c^{\mathrm{far}}|
 \le C\sum_{i\ge0}2^{-3i}\kappa^{-3}\frac R{K^3}\norm{g_K}_2
 \le C\frac R{K^3}\norm{g_K}_2.
\end{equation}
Together,
\begin{equation}\label{eq:highleak}
 \boxed{\quad\sup_{x\in\Gamma_{c,R}}|q_c(x)|^2
 \le C\frac{R^2}{K^6}\norm{g_K}_2^2.\quad}
\end{equation}
Unlike \eqref{eq:lowleak}, this is a global input-norm estimate; it does not require summing local data over individual high-frequency sectors.

By \eqref{eq:AH}, the bounded overlap of $\Gamma_{c,R}$, and the multiplicity bound $\sum_c\norm{f_c}_2^2\le CM^2\norm{g_L}_2^2$,
\begin{equation}\label{eq:lowmass}
 \sum_c\int_{\Gamma_{c,R}}|Ef_c|^2dx\le CM^2R\norm{g_L}_2^2.
\end{equation}
The same bound holds with $|Ef_c|^2$ replaced by $|\ell_c|^2$, using
$|\ell_c|^2\le2|Eg_L|^2+2|Ef_c|^2$, bounded sector overlap, and \eqref{eq:AH}. Hence \eqref{eq:highleak} gives
\begin{equation}\label{eq:IIIIV}
 \sum_c\int_{\Gamma_{c,R}}(III_{c,R}+IV_{c,R})\,dx
 \le CM^2\frac{R^3}{K^6}\norm{g_L}_2^2\norm{g_K}_2^2.
\end{equation}
Summing $R<L^2$ yields
\begin{equation}\label{eq:IIIIVall}
 \sum_{R,c}\int_{\Gamma_{c,R}}(III_{c,R}+IV_{c,R})\,dx
 \le CM^2\left(\frac LK\right)^6\norm{g_L}_2^2\norm{g_K}_2^2.
\end{equation}
The fourth, simultaneous-leakage term is therefore included. There is no independent turning-point estimate to prove: the high-frequency turning-point directions are retained in $h_c$, where only Lemmas \ref{lem:fiber} and \ref{lem:tube}, valid for arbitrary $L^2$ data on their stated supports, are used.

\subsection{Assembly}
Combining \eqref{eq:inner}, \eqref{eq:far}, \eqref{eq:mainall}, \eqref{eq:IIall}, and \eqref{eq:IIIIVall}, with \eqref{eq:four}, proves
\begin{equation}\label{eq:bilinearfinal}
 \norm{Eg_LEg_K}_2^2\le
 C\left(\kappa M^3\log(2\kappa)+\kappa M^{-4}+M^2+1\right)
 \frac LK\norm{g_L}_2^2\norm{g_K}_2^2.
\end{equation}
Since $\kappa,M$ were fixed absolutely in \eqref{eq:constants}, this is Theorem \ref{thm:bilinear} (the upper bound). The deduction \eqref{eq:LP4}--\eqref{eq:discrete} completes the proof of Theorem \ref{thm:main}.

\section{Sharpness of the band separation}
\subsection{Why the band-separation factor is optimal}
We show that the factor $L/K$ in \eqref{eq:bilinear} cannot be
replaced by a uniformly smaller order. The calculation also explains
why the numerator is the lower angular degree rather than an artifact
of the proof. All the data below are normalized in $L^2(\Sph)$.

A preliminary example isolates the decay in $K$. Take the constant
harmonic $g_1=Y_{00}$ and any normalized $g_K\in\HH_K$. Since
$Eg_1(r\theta)=\sqrt{4\pi}\sin r/r$ and
$Eg_K(r\theta)=4\pi i^Kj_K(r)g_K(\theta)$, polar integration gives
\begin{equation}\label{eq:constantlower}
 \norm{Eg_1Eg_K}_2^2
 =C\int_0^\infty\sin^2r\,|j_K(r)|^2\,dr
 \gtrsim K^{-1}.
\end{equation}
We explain the radial lower bound below. This example shows that
$K^{-1}$ cannot be improved when the low frequency is fixed, but does
not yet exhibit the factor $L$.

\smallskip\noindent\emph{Construction of concentrated low-band data.}
Let $N$ be the north pole and let $L\to\infty$ through integers. For
each even $k\in[L,2L-1]$, put
\[
 Z_k(N,\omega)=\frac{2k+1}{4\pi}P_k(N\cdot\omega),
\]
where $P_k$ is the Legendre polynomial. Define
\begin{equation}\label{eq:zonal}
 H_L=\sum_{\substack{L\le k\le2L-1\\k\text{ even}}}Z_k(N,\cdot),
 \qquad g_L=H_L/\norm{H_L}_2.
\end{equation}
The kernels in different degrees are orthogonal, and
$\norm{Z_k(N,\cdot)}_2^2=Z_k(N,N)=(2k+1)/(4\pi)$.
There are comparable to $L$ even degrees in the sum; hence
\begin{equation}\label{eq:zonalnorm-sharp}
 \norm{H_L}_2^2
 =\sum_{\substack{L\le k\le2L-1\\k\text{ even}}}
   \frac{2k+1}{4\pi}\asymp L^2.
\end{equation}
For $0\le\theta\le\delta/L$, the elementary bound
$|P_k'(t)|\le k(k+1)/2$ and $1-\cos\theta\le\theta^2/2$
give $P_k(\cos\theta)\ge1-C\delta^2$ for all degrees in the
sum. Choose $\delta>0$ fixed and small so that this is at least
$1/2$. Every summand is then positive. It follows that
\begin{equation}\label{eq:zonallocal}
 |g_L(\omega)|\ge cL,
 \qquad \omega\in\mathcal C_L
 :=\{\omega:d(\omega,N)\le\delta/L\}.
\end{equation}

\smallskip\noindent\emph{High-degree mass in the same cap.}
Let $g_K=Y_{K0}$ be the $L^2$-normalized zonal harmonic with pole
$N$, where $K/L^2\to\infty$. We claim
\begin{equation}\label{eq:Kmass}
 m_K(L):=\int_{\mathcal C_L}|Y_{K0}(\omega)|^2\,d\sigma(\omega)
 \ge cL^{-1}.
\end{equation}
To see the scale of the integral, restrict to the subannulus
$\delta/(2L)\le\theta\le\delta/L$. Here
$(K+1/2)\theta\to\infty$ uniformly. The uniform Hilb--Bessel
asymptotic for Legendre functions, followed by the large-argument
asymptotic of $J_0$ (see \cite{dlmf}), yields on that annulus
\[
 P_K(\cos\theta)
 =\left(\frac{2}{\pi(K+\tfrac12)\sin\theta}\right)^{1/2}
  \left[\cos\bigl((K+\tfrac12)\theta-\tfrac\pi4\bigr)
      +o(1)\right],
\]
where the error is uniform as $K/L\to\infty$. Since
$|Y_{K0}(\theta)|^2=(2K+1)|P_K(\cos\theta)|^2/(4\pi)$,
and since the area element is $2\pi\sin\theta\,d\theta$,
we obtain
\[
 \begin{aligned}
 m_K(L)&\ge c\int_{\delta/(2L)}^{\delta/L}
 \Bigl[\cos^2\bigl((K+\tfrac12)\theta-\tfrac\pi4\bigr)
              -o(1)\Bigr]d\theta\\
 &\ge c\left(\frac{\delta}{4L}-\frac{C}{K}
                     -\frac{o(1)}{L}\right)
 \ge \frac{c'}L.
 \end{aligned}
\]
For the middle step we used
$\int_a^b\cos^2(t\theta+\alpha)\,d\theta
=(b-a)/2+O(t^{-1})$. The constant $c'$ is independent of $L,K$
when $K/L^2$ is sufficiently large. Combining this with
\eqref{eq:zonallocal} gives
\begin{equation}\label{eq:zonal-product}
 \int_{\mathcal C_L}|g_L|^2|g_K|^2\,d\sigma
 \ge cL^2m_K(L)\gtrsim L.
\end{equation}

\smallskip\noindent\emph{Uniform outgoing expansion of the low band.}
The exact finite spherical-Hankel representation
\eqref{eq:hankel}, and the fact that every degree of $g_L$ is even,
give
\begin{equation}\label{eq:sharpasympt}
 rEg_L(r\theta)=4\pi\sin r\,g_L(\theta)
                     +\mathcal R_L(r,\theta).
\end{equation}
In order to use this identity on the small cap we record a
\emph{pointwise} error bound. The addition theorem on the entire
band, together with $\norm{g_L}_2=1$, shows that a diagonal degree
multiplier $u_n$ satisfies
\[
 \norm{u_n}_\infty
 \le C L\max_{k\in[L,2L-1]}a_n(k)
 \le CL\frac{(CL^2)^n}{n!}.
\]
The same estimate holds for the antipodal contribution in the Hankel
formula. Summing $n\ge1$ gives
\begin{equation}\label{eq:sharp-error-pt}
 \sup_{\theta\in\Sph}|\mathcal R_L(r,\theta)|
 \le CL\bigl(e^{CL^2/r}-1\bigr).
\end{equation}
Fix a sufficiently large constant $A$ and suppose $K\ge AL^2$.
For $10K\le r\le20K$ the right-hand side of
\eqref{eq:sharp-error-pt} is at most $\varepsilon L$, where
$\varepsilon$ can be made arbitrarily small by increasing $A$.
Using $|a+b|^2\ge\tfrac12|a|^2-|b|^2$, we therefore get on
$\mathcal C_L$
\begin{equation}\label{eq:sharp-pointwise-square}
 |rEg_L(r\theta)|^2
 \ge cL^2\sin^2r-\varepsilon^2L^2.
\end{equation}

\smallskip\noindent\emph{Radial integration.}
We need two standard uniform consequences of the spherical-Bessel
oscillatory asymptotic in the region $10K\le r\le20K$:
\begin{equation}\label{eq:sharpbessel}
 \int_{10K}^{20K}|j_K(r)|^2\,dr\le \frac C K,
 \qquad
 \int_{10K}^{20K}\sin^2r\,|j_K(r)|^2\,dr\ge\frac c K.
\end{equation}
Here is why the second estimate is uniform, despite the two
oscillatory factors. On this interval the Debye formula has the form
\[
 j_K(r)=r^{-1}a_K(r)
     \bigl[\cos\Phi_K(r)+O(K^{-1})\bigr],
 \qquad c\le a_K(r)\le C,
\]
with $\Phi_K'(r)$ bounded away from zero and with
$1+\Phi_K'(r)$ bounded away from zero. The first assertion of
\eqref{eq:sharpbessel} follows at once. For the second, expand
\[
 \sin^2r\cos^2\Phi_K(r)
 =\frac14\bigl(1-\cos2r+\cos2\Phi_K
               -\cos2r\cos2\Phi_K\bigr).
\]
The integrals of $\cos2r$, $\cos2\Phi_K$, and
$\cos2(r+\Phi_K)$ are $O(1)$ by elementary integration by parts
and the uniform Debye derivative bounds. The remaining term
$\cos2(r-\Phi_K)$ is bounded below by $-1$ pointwise.
Since
$\cos2r\cos2\Phi_K
=\tfrac12\cos2(r-\Phi_K)+\tfrac12\cos2(r+\Phi_K)$,
this yields
\[
 \int_{10K}^{20K}\sin^2r\cos^2\Phi_K(r)\,dr
 \ge cK-O(1).
\]
The $O(K^{-1})$ asymptotic error does not affect the lower bound
for large $K$; using $r\asymp K$ proves \eqref{eq:sharpbessel}.
The same argument proves \eqref{eq:constantlower}.

Finally, $Eg_K(r\theta)=4\pi i^Kj_K(r)g_K(\theta)$.
Multiplying \eqref{eq:sharp-pointwise-square} by
$|j_K(r)|^2|g_K(\theta)|^2$, integrating over
$[10K,20K]\times\mathcal C_L$, and using
\eqref{eq:sharpbessel} gives, provided $\varepsilon$ was chosen
small enough,
\[
 \begin{aligned}
 \norm{Eg_LEg_K}_2^2
 &\ge (4\pi)^2\int_{10K}^{20K}|j_K(r)|^2
       \int_{\mathcal C_L}|rEg_L(r\theta)|^2
                       |g_K(\theta)|^2\,d\sigma(\theta)\,dr\\
 &\ge cL^2m_K(L)
      \left[\int_{10K}^{20K}\sin^2r|j_K(r)|^2\,dr
       -C\varepsilon^2\int_{10K}^{20K}|j_K(r)|^2\,dr\right]\\
 &\gtrsim L^2\cdot L^{-1}\cdot K^{-1}=\frac LK.
 \end{aligned}
\]
This proves the sharpness assertion in Theorem~\ref{thm:bilinear}.
Notice that the construction lies in the far field and does not rely
on any of the Fresnel localization estimates.

\section{Sharp multiscale data and optimality of the Besov exponent}
\label{sec:multiscale}
The factor $L/K$ is optimal by the construction in Section~6. We now
show that the fourth-power summation in the main theorem is optimal
in a stronger sense: there are many-band configurations on which it
computes the spatial $L^4$ norm up to two fixed constants. This is the
appropriate example for distinguishing our refinement from an
application of the triangle inequality to individual bands.

\begin{lemma}\label{lem:band-sharp}
There are $L_0\ge1$ and $c_0>0$ such that, for each dyadic
$L\ge L_0$, one can find $h_L\in L^2(\Sph)$ with
\begin{equation}\label{eq:band-sharp}
 \norm{h_L}_2=1,\qquad
 \supp_\Lambda h_L\subset[L/2,2L],\qquad
 \norm{Eh_L}_{L^4(\RR^3)}\ge c_0.
\end{equation}
The functions may be chosen to involve only even spherical-harmonic degrees.
\end{lemma}
\begin{proof}
Fix a pole $N\in\Sph$ and, for even integers $k\in[L,2L-1]$,
write the zonal reproducing kernel
\[
 Z_k(N,\omega)=\frac{2k+1}{4\pi}P_k(N\cdot\omega).
\]
Set $H_L=\sum_{L\le k\le2L-1,\ k\text{ even}}Z_k(N,\cdot)$ and
$h_L=H_L/\norm{H_L}_2$. Orthogonality and the addition theorem give
\begin{equation}\label{eq:hLnorm}
 \norm{H_L}_2^2=H_L(N)
 =\sum_{\substack{L\le k\le2L-1\\k\text{ even}}}
   \frac{2k+1}{4\pi}\asymp L^2.
\end{equation}
There are therefore fixed constants $a,b>0$ such that
$aL\le h_L(N)\le bL$. Moreover, the elementary estimate
$|P_k'(t)|\le k(k+1)/2$ for $-1\le t\le1$ implies
$P_k(\cos\theta)\ge1-Ck^2\theta^2$. Choose a sufficiently small
constant $\delta>0$. For $d(\omega,N)\le\delta/L$, every Legendre
factor in the sum is at least $1/2$, so
$|h_L(\omega)|\ge cL$ on that cap. Its area is comparable to $L^{-2}$;
consequently
\begin{equation}\label{eq:hL4}
 \norm{h_L}_{L^4(\Sph)}\ge cL^{1/2}.
\end{equation}
The reproducing-kernel Cauchy--Schwarz estimate also gives
$\norm{h_L}_\infty\le CL$. Interpolation between $L^2$ and
$L^\infty$ therefore yields the matching upper bound
$\norm{h_L}_4\le CL^{1/2}$.

We turn to the extension. Because all harmonic degrees of $h_L$
are even, the first term of the exact finite Hankel expansion
\eqref{eq:hankel} takes the form
\begin{equation}\label{eq:farband-leading}
 rEh_L(r\theta)=4\pi\sin r\,h_L(\theta)+R_L(r,\theta).
\end{equation}
For the $n$th remainder multiplier $u_n$ in that expansion,
its spectral coefficients satisfy
$a_n(k)\le(C L^2)^n/n!$ on the support of $h_L$. Since $u_n$
has degrees at most $2L$, the addition theorem and $L^2$--$L^\infty$
interpolation imply
\[
 \norm{u_n}_{L^4(\Sph)}
 \le CL^{1/2}\norm{u_n}_2
 \le CL^{1/2}\frac{(CL^2)^n}{n!}.
\]
Summing the nonzero terms of the finite expansion and using the same
estimate for the antipodal terms gives, for $r\ge AL^2$,
\begin{equation}\label{eq:farband-error}
 \norm{R_L(r,\cdot)}_{L^4(\Sph)}
 \le CL^{1/2}\big(e^{C/A}-1\big).
\end{equation}
Choose the fixed number $A$ large enough that the right side is at
most $\pi cL^{1/2}$, where $c$ is the constant in
\eqref{eq:hL4}. Whenever $|\sin r|\ge1/2$ and
$AL^2\le r\le2AL^2$, the triangle inequality in $L^4(\Sph)$
then gives
\[
 \norm{rEh_L(r\cdot)}_{L^4(\Sph)}\ge c'L^{1/2}.
\]
For all sufficiently large $L$, the set of such radii has measure
at least a fixed proportion of $AL^2$. Therefore
\[
 \begin{aligned}
 \norm{Eh_L}_{L^4(\RR^3)}^4
 &\ge \int_{AL^2}^{2AL^2}r^{-2}
        \norm{rEh_L(r\cdot)}_{L^4(\Sph)}^4\,dr\\
 &\gtrsim L^2\int_{AL^2}^{2AL^2}
             r^{-2}\one_{\{|\sin r|\ge1/2\}}\,dr
 \gtrsim1.
 \end{aligned}
\]
This proves the claim. Notice that the argument uses the exact
spherical-Hankel expansion already needed for the main theorem;
no additional paraboloid limit is required.
\end{proof}

\begin{theorem}\label{thm:exact-ell4}
There are smooth, normalized spherical data $h_m$ with mutually
disjoint angular spectral bands and constants $0<c<C<\infty$ such
that, for every finite scalar sequence $(a_m)$,
\begin{equation}\label{eq:exact-ell4}
 c\left(\sum_m|a_m|^4\right)^{1/4}
 \le \norm{E\left(\sum_ma_mh_m\right)}_{L^4(\RR^3)}
 \le C\left(\sum_m|a_m|^4\right)^{1/4}.
\end{equation}
Here $B^0_{2,q}$ denotes the analogous inhomogeneous Besov space
with $\ell^q$ summation over dyadic angular blocks. On the same
family, for each $1\le q\le\infty$,
\begin{equation}\label{eq:besov-ellq}
 \norm{\sum_m a_mh_m}_{B^0_{2,q}(\Sph)}
 \asymp\norm{(a_m)}_{\ell^q}.
\end{equation}
In particular, no bounded extension
$E:B^0_{2,q}(\Sph)\to L^4(\RR^3)$ exists for $q>4$.
\end{theorem}
\begin{proof}
Choose the functions in Lemma~\ref{lem:band-sharp} at scales
$L_m=16^{m+m_0}$, with $m_0$ large enough that the lemma applies.
Their spectra lie in the disjoint intervals $[L_m,2L_m]$.
We first verify \eqref{eq:besov-ellq} explicitly, including the
uniformity in the number of bands. Every smooth dyadic multiplier
$P_j$ is supported in an angular-frequency interval of fixed
multiplicative width. Since $L_{m+1}/L_m=16$, no such interval
intersects two distinct chosen bands. For each $m$, let $J_m$ be
the set of indices $j$ for which $P_jh_m$ is nonzero. The
cardinality of $J_m$ is bounded by an absolute constant, say $N_0$,
and the sets $J_m$ are disjoint. The $L^2$ spectral multiplier bound
gives $\norm{P_jh_m}_2\le C$ for $j\in J_m$. Conversely,
$h_m=\sum_{j\in J_m}P_jh_m$ because the inhomogeneous multipliers
form a partition of unity. Thus
\[
 1=\norm{h_m}_2\le\sum_{j\in J_m}\norm{P_jh_m}_2,
\]
so at least one of these finitely many block norms is at least
$N_0^{-1}$. For $g=\sum_ma_mh_m$ we consequently have
\[
 c|a_m|\le
 \Big(\sum_{j\in J_m}\norm{P_jg}_2^q\Big)^{1/q}
 \le C|a_m|,
\]
for $1\le q<\infty$, with the analogous maximum inequality when
$q=\infty$. Summing over the disjoint $J_m$ proves
\eqref{eq:besov-ellq}, with constants independent of how many
coefficients are nonzero.

For the upper bound in \eqref{eq:exact-ell4}, apply
Theorem~\ref{thm:main} and \eqref{eq:besov-ellq} with $q=4$.
For the lower bound, choose smooth multipliers
$Q_m=\chi(\Lambda/L_m)$ which equal one on $[L_m,2L_m]$
and vanish on the other selected bands. Thus, for
$g=\sum a_mh_m$,
$Q_mg=a_mh_m$, and, at every radius,
$Q_m(Eg(r\cdot))=a_mEh_m(r\cdot)$.
The reverse spherical Littlewood--Paley square-function bound,
valid uniformly for this lacunary family, yields
\[
 \begin{aligned}
 \norm{Eg}_{L^4(\RR^3)}^4
 &\ge c\int_{\RR^3}
       \left(\sum_m|Q_m(Eg)(x)|^2\right)^2dx\\
 &\ge c\sum_m|a_m|^4\norm{Eh_m}_{L^4(\RR^3)}^4
 \ge c'\sum_m|a_m|^4.
 \end{aligned}
\]
For clarity, the square-function inequality used in the first step
has a direct proof in this lacunary setting. For each choice of signs
$\varepsilon_m\in\{-1,1\}$, the symbol
$\sum_m\varepsilon_m\chi(\lambda/L_m)$ is a smooth lacunary
spherical multiplier with bounds uniform in the sign choices and the
number of terms. The spherical multiplier theorem of
Strichartz and Bonami--Clerc therefore gives, at every fixed radius,
\[
 \Big\|\sum_m\varepsilon_m Q_mF\Big\|_{L^4(\Sph)}
 \le C\|F\|_{L^4(\Sph)}.
\]
Averaging the fourth power over independent random signs and applying
the scalar Khintchine inequality yields
\[
 \left\|\Big(\sum_m|Q_mF|^2\Big)^{1/2}\right\|_{L^4(\Sph)}
 \le C\|F\|_{L^4(\Sph)}.
\]
Integration in the radius supplies the displayed reverse bound above.
The next step uses only nonnegativity of the mixed terms, and
Lemma~\ref{lem:band-sharp} supplies the final uniform lower bound.

Taking $a_m=1$ on $1\le m\le J$ gives
$\norm{Eg}_4\asymp J^{1/4}$, whereas
$\norm{g}_{B^0_{2,q}}\asymp J^{1/q}$ for finite $q$
(and $\asymp1$ for $q=\infty$).
The ratio diverges when $q>4$, proving the final assertion.
\end{proof}

For this construction the older $L^2$ Stein--Tomas bound has a
right-hand side of order $J^{1/2}$, while the true $L^4$ size is
$J^{1/4}$. This discrepancy grows without bound. The example also
cannot be controlled by summing the individual extension norms,
since that would give order $J$. Thus the $\ell^4$ improvement is
a genuinely multiscale effect, not a repackaging of a one-band estimate.

\section{Cubic nonlinear Helmholtz waves and scattering}\label{sec:applications}
The refined restriction estimate changes the smallness condition in
an established nonlinear fixed-point problem. We work with the outgoing
fundamental solution
\begin{equation}\label{eq:out-resolvent}
 \Phi_+(x)=\frac{e^{i|x|}}{4\pi|x|},\qquad
 \mathcal R_+F=\Phi_+*F,
\end{equation}
where the convolution is understood first for test functions and then
by continuous extension. The limiting-absorption estimate
\begin{equation}\label{eq:L43L4-resolvent}
 \norm{\mathcal R_+F}_{L^4(\RR^3)}
 \le C_R\norm{F}_{L^{4/3}(\RR^3)}
\end{equation}
is the standard resolvent input for the cubic fixed-point argument
\cite{gutierrez2004}. We do not claim this resolvent bound as a new result.
The novel step is that the incident field $Eg$ can now be controlled
in terms of a strictly weaker angular input norm.

\subsection{Small solutions from angular Besov data}
\begin{theorem}\label{thm:nonlinear}
Let $Q\in L^\infty(\RR^3;\mathbb C)$ and write $Q_0=\norm Q_\infty$.
There is $\varepsilon_0>0$, depending on $C_R$, the extension constant,
and $Q_0$, such that whenever
$g\in B^0_{2,4}(\Sph)$ and
$\norm g_{B^0_{2,4}}\le\varepsilon_0$, there is a unique
solution in a sufficiently small ball of $L^4(\RR^3)$ to
\begin{equation}\label{eq:nl-fixed-point}
 u=Eg+\mathcal R_+\bigl(Q|u|^2u\bigr).
\end{equation}
It satisfies $(-\Delta-1)u=Q|u|^2u$ in distributions and
\begin{equation}\label{eq:nl-size}
 \norm u_4\le C\norm g_{B^0_{2,4}},\qquad
 \norm{u-Eg}_4\le C Q_0\norm g_{B^0_{2,4}}^3.
\end{equation}
The solution map $g\mapsto u$ is locally Lipschitz and
real-analytic near the origin, when the complex Banach spaces are
regarded as real spaces.
\end{theorem}
\begin{proof}
Let $C_E$ be the constant in \eqref{eq:main}, put
$v=Eg$ and $\delta=C_E\norm g_{B^0_{2,4}}$, and write
$b=C_RQ_0$. Theorem~\ref{thm:main} implies $\norm v_4\le\delta$.
For $u\in L^4$ define
\[
 \mathcal T_v(u)=v+\mathcal R_+\bigl(Q|u|^2u\bigr).
\]
The map is well defined on $L^4$: by H\"older's inequality,
$\norm{|u|^2u}_{4/3}=\norm u_4^3$, and the resolvent estimate
\eqref{eq:L43L4-resolvent} therefore gives
\begin{equation}\label{eq:nl-contraction1}
 \norm{\mathcal T_v(u)}_4\le\delta+b\norm u_4^3.
\end{equation}
The pointwise algebraic inequality
\[
 \bigl||z|^2z-|w|^2w\bigr|
 \le 3(|z|^2+|w|^2)|z-w|,
 \qquad z,w\in\mathbb C,
\]
followed by H\"older's inequality with exponents $2$ and $4$,
shows that
\begin{equation}\label{eq:nl-cubic-difference}
 \norm{|u|^2u-|w|^2w}_{4/3}
 \le 3(\norm u_4^2+\norm w_4^2)\norm{u-w}_4.
\end{equation}
Consequently,
\begin{equation}\label{eq:nl-contraction2}
 \norm{\mathcal T_v(u)-\mathcal T_v(w)}_4
 \le3b(\norm u_4^2+\norm w_4^2)\norm{u-w}_4.
\end{equation}

Consider the closed ball $\mathbb B_{2\delta}$ of radius $2\delta$
in $L^4$. For $u\in\mathbb B_{2\delta}$, estimate
\eqref{eq:nl-contraction1} gives
\[
 \norm{\mathcal T_v(u)}_4\le
 \delta+8b\delta^3.
\]
Thus $\mathcal T_v$ maps this ball into itself whenever
$8b\delta^2\le1$. For $u,w$ in the same ball,
\eqref{eq:nl-contraction2} gives the Lipschitz constant
$24b\delta^2$. Choose $\varepsilon_0>0$ so that
$\delta\le C_E\varepsilon_0$ and both
$8b\delta^2\le1$ and $24b\delta^2\le\tfrac12$ hold.
(When $Q=0$ the equation is linear and these conditions are
unnecessary.) The contraction mapping theorem now yields a unique
fixed point in $\mathbb B_{2\delta}$. The word \emph{unique} is
restricted to this small ball; no assertion about all possibly large
$L^4$ solutions is needed.

The first estimate in \eqref{eq:nl-size} follows from
$\norm u_4\le2\delta$. The fixed-point equation and the resolvent
bound give the quantitative correction estimate
\[
 \norm{u-v}_4
 =\norm{\mathcal R_+(Q|u|^2u)}_4
 \le b\norm u_4^3
 \le8b\delta^3,
\]
which proves the second inequality in \eqref{eq:nl-size}.
Since $(-\Delta-1)\Phi_+=\delta_0$ distributionally,
applying $(-\Delta-1)$ to the fixed-point equation yields
$(-\Delta-1)u=Q|u|^2u$ in the sense of distributions; the incident
field satisfies the homogeneous Helmholtz equation term by term.

We next make the dependence on $g$ precise. Let $g,h$ be two data
in the same sufficiently small Besov ball, and denote by $u_g,u_h$
their fixed points. Subtracting their equations yields
\[
 \begin{aligned}
 \norm{u_g-u_h}_4
 &\le\norm{E(g-h)}_4
 +3b(\norm{u_g}_4^2+\norm{u_h}_4^2)
                \norm{u_g-u_h}_4.
 \end{aligned}
\]
The coefficient of the final term is at most $1/2$ by our choice of
ball. Move it to the left to get
\[
 \norm{u_g-u_h}_4
 \le2\norm{E(g-h)}_4
 \le2C_E\norm{g-h}_{B^0_{2,4}}.
\]
Finally, regard the complex Banach spaces as real ones and define
\[
 \mathcal F(g,u)=u-Eg-\mathcal R_+(Q|u|^2u).
\]
The map $u\mapsto |u|^2u$ is a continuous real cubic polynomial
from $L^4$ to $L^{4/3}$, so $\mathcal F$ is real-analytic near
$(0,0)$. Its derivative with respect to $u$ at $(0,0)$ is the
identity on $L^4$. The real-analytic implicit-function theorem
therefore produces a real-analytic local branch $g\mapsto u_g$.
By uniqueness in the small ball, this is precisely the branch
constructed above.
\end{proof}

\subsection{The scattering amplitude has improved angular summability}
For $F\in L^{4/3}(\RR^3)$, the dual restriction theorem makes
$\widehat F|_{\Sph}$ a well-defined element of
$B^0_{2,4/3}(\Sph)$, even when $F$ is not integrable.
We therefore define the on-shell scattered amplitude of the solution
in Theorem~\ref{thm:nonlinear} by
\begin{equation}\label{eq:on-shell}
 \mathcal A(g)=
 \left.\widehat{Q|u_g|^2u_g}\right|_{\Sph}.
\end{equation}
For compactly supported $Q$, this is also the classical outgoing
far-field coefficient, with the normalization described below.

\begin{theorem}\label{thm:amplitude}
Under the hypotheses of Theorem~\ref{thm:nonlinear},
\begin{equation}\label{eq:amplitude-besov}
 \boxed{\quad
 \norm{\mathcal A(g)}_{B^0_{2,4/3}(\Sph)}
 \le C Q_0\norm g_{B^0_{2,4}(\Sph)}^3.
 \quad}
\end{equation}
For sufficiently small $g,h$ in the same Besov ball,
\begin{equation}\label{eq:amplitude-lipschitz}
 \norm{\mathcal A(g)-\mathcal A(h)}_{B^0_{2,4/3}}
 \le C Q_0\bigl(\norm g_{B^0_{2,4}}+\norm h_{B^0_{2,4}}\bigr)^2
             \norm{g-h}_{B^0_{2,4}}.
\end{equation}
The map $\mathcal A$ is real-analytic near the origin.
If $Q$ has bounded support, the scattered field satisfies
\begin{equation}\label{eq:amplitude-farfield}
 (u_g-Eg)(r\theta)
 =\frac{e^{ir}}{4\pi r}\mathcal A(g)(\theta)+O_g(r^{-2})
\end{equation}
uniformly for $\theta\in\Sph$ as $r\to\infty$.
\end{theorem}
\begin{proof}
By Theorem~\ref{thm:nonlinear} and H\"older's inequality,
\[
 \norm{Q|u_g|^2u_g}_{L^{4/3}}
 \le Q_0\norm{u_g}_4^3
 \le C Q_0\norm g_{B^0_{2,4}}^3.
\]
The dual restriction estimate \eqref{eq:dual} applied to this
source yields \eqref{eq:amplitude-besov}.
For the difference estimate, use
$\norm{|u|^2u-|w|^2w}_{4/3}
\le C(\norm u_4+\norm w_4)^2\norm{u-w}_4$
and the local Lipschitz bound from Theorem~\ref{thm:nonlinear}.
Real-analyticity follows by composition with the real cubic map and
the bounded linear restriction operator.

If $Q$ is supported in a ball $B(0,R_Q)$, then
$F=Q|u_g|^2u_g$ belongs to $L^1$ as well as $L^{4/3}$,
because H\"older gives
$\norm F_1\le |B(0,R_Q)|^{1/4}\norm F_{4/3}$.
The convolution formula \eqref{eq:out-resolvent} is therefore
absolutely convergent for all sufficiently large $|x|$.
For $x=r\theta$ with $r>2R_Q$, expand
\[
 |r\theta-y|=r-\theta\cdot y+O_{R_Q}(r^{-1}),\qquad
 |r\theta-y|^{-1}=r^{-1}+O_{R_Q}(r^{-2})
\]
uniformly for $|y|\le R_Q$. Substitution into the kernel gives
\[
 \mathcal R_+F(r\theta)
 =\frac{e^{ir}}{4\pi r}
    \int e^{-i\theta\cdot y}F(y)\,dy+O_{F,R_Q}(r^{-2}),
\]
which is \eqref{eq:amplitude-farfield}.
For arbitrary bounded $Q$, we claim only the on-shell Besov
statement, not a classical pointwise far-field asymptotic.
\end{proof}

\subsection{Incident data with unbounded square-integrable energy}
The essential improvement is quantitative and multiscale, not
merely the existence of data having large $L^2$ norm and small
extension norm. Indeed, a single high-degree harmonic can have
small $L^4$ extension relative to its $L^2$ norm. Our example instead
uses one near-extremal band at each of many separated scales, so that
the \emph{two-sided} estimate \eqref{eq:exact-ell4} identifies the
actual incident-wave size. The following statement records the
resulting enlargement of the nonlinear small-data regime.

\begin{corollary}\label{cor:large-L2}
For every $M>0$ and every $\eta>0$ sufficiently small, there is a
smooth angular density $g\in L^2(\Sph)$ such that
\[
 \norm{g}_{L^2(\Sph)}>M,
 \qquad \norm g_{B^0_{2,4}(\Sph)}<\eta,
\]
and whose nonlinear solution from Theorem~\ref{thm:nonlinear}
satisfies $\norm u_4\le C\eta$.
There are also distributions $g\in B^0_{2,4}(\Sph)\setminus
L^2(\Sph)$ of arbitrarily small Besov norm for which the same
small solution exists.
\end{corollary}
\begin{proof}
Let $h_m$ be the normalized, pairwise orthogonal data from
Theorem~\ref{thm:exact-ell4}. For finite $J$ and scalar
$\varepsilon>0$, put
\[
 g_{\varepsilon,J}=\varepsilon\sum_{m=1}^{J}h_m.
\]
The disjoint spectral supports give
$\norm{g_{\varepsilon,J}}_2=\varepsilon J^{1/2}$,
whereas \eqref{eq:besov-ellq} gives
$\norm{g_{\varepsilon,J}}_{B^0_{2,4}}\asymp
\varepsilon J^{1/4}$, uniformly in $J$.
Taking, for example, $J=\lfloor\varepsilon^{-3}\rfloor$
and sending $\varepsilon\downarrow0$ gives
\[
 \norm{g_{\varepsilon,J}}_{B^0_{2,4}}\asymp
 \varepsilon^{1/4}\longrightarrow0,
 \qquad
 \norm{g_{\varepsilon,J}}_2\asymp
 \varepsilon^{-1/2}\longrightarrow\infty.
\]
This proves the first assertion and the claim about the solution
by Theorem~\ref{thm:nonlinear}. More quantitatively, the two-sided
law \eqref{eq:exact-ell4} gives
$\|Eg_{\varepsilon,J}\|_4\asymp\varepsilon J^{1/4}$,
while the classical Stein--Tomas estimate bounds this expression
only by $C\varepsilon J^{1/2}$; summing the one-band extension
bounds separately gives merely $C\varepsilon J$. The discrepancy
between the correct size and the $L^2$ estimate is an unbounded
factor $J^{1/4}$. For infinite-energy data, let
\[
 g=\varepsilon\sum_{m=1}^{\infty}m^{-1/3}h_m.
\]
The series converges in $B^0_{2,4}$, because
$\sum_m m^{-4/3}<\infty$, but its terms are orthogonal in $L^2$
and $\sum_m m^{-2/3}=\infty$. Consequently
$g\in B^0_{2,4}\setminus L^2$. Taking $\varepsilon$ small
ensures that Theorem~\ref{thm:nonlinear} applies. For this series,
$\sum_m m^{-1/3}\|Eh_m\|_4$ diverges by the one-band
lower bound. Thus neither Stein--Tomas, which requires $L^2$ data,
nor the termwise triangle estimate establishes the global
$L^4$ conclusion; the refined multiscale estimate does.
\end{proof}

\subsection{Decoupling across separated angular scales}
The separated-band estimate controls cross-scale nonlinear
interactions at the first nonzero order of the scattering map.
It is important that the resulting gain is asserted for the leading
cubic term; the higher-order corrections do not generally decay
with angular separation.

\begin{theorem}\label{thm:scattering-decouple}
Suppose $g_L,g_K$ have spectra as in \eqref{eq:bands},
$K\ge C_0L$, and
$\norm{g_L}_2+\norm{g_K}_2\le\varepsilon$,
where $\varepsilon$ is sufficiently small for Theorem~\ref{thm:nonlinear}.
Then
\begin{equation}\label{eq:scattering-decouple}
 \boxed{\quad
 \norm{\mathcal A(g_L+g_K)-\mathcal A(g_L)-\mathcal A(g_K)}_{B^0_{2,4/3}}
 \le C_Q\left[\left(\frac LK\right)^{1/2}
 \varepsilon^3+\varepsilon^5\right].\quad}
\end{equation}
Here $C_Q$ depends on the fixed coefficient bound $\norm Q_\infty$
and the linear extension and resolvent constants, not on $L,K$.
In particular the leading cubic part of the scattering map is
asymptotically additive as $K/L\to\infty$.
\end{theorem}
\begin{proof}
Put $v=Eg_L$, $w=Eg_K$, and define the cubic Born approximation
\[
 \mathcal A_3(g)=
 \left.\widehat{Q|Eg|^2Eg}\right|_{\Sph}.
\]
We first compare the exact scattering amplitude with its Born
approximation, and then estimate the cross terms in the cubic
expression. This separates the small-data expansion from the angular
scale separation.

\smallskip\noindent\emph{The fifth-order error.}
For each of the three inputs $g=g_L,g_K,g_L+g_K$, the extension
estimate and the triangle inequality give
$\norm{Eg}_4\le C\varepsilon$. Theorem~\ref{thm:nonlinear}
therefore yields
\[
 \norm{u_g}_4\le C\varepsilon,
 \qquad\norm{u_g-Eg}_4\le C_Q\varepsilon^3.
\]
Apply the cubic difference estimate \eqref{eq:nl-cubic-difference}
to $u_g$ and $Eg$, then use the dual restriction estimate:
\[
 \begin{aligned}
 \norm{\mathcal A(g)-\mathcal A_3(g)}_{B^0_{2,4/3}}
 &\le C Q_0
   \norm{|u_g|^2u_g-|Eg|^2Eg}_{4/3}\\
 &\le C Q_0(\norm{u_g}_4^2+\norm{Eg}_4^2)
             \norm{u_g-Eg}_4\\
 &\le C_Q\varepsilon^5.
 \end{aligned}
\]
The constants here are uniform in $L,K$. In particular we obtain
\begin{equation}\label{eq:Born-error}
 \norm{\mathcal A(g)-\mathcal A_3(g)}_{B^0_{2,4/3}}
 \le C_Q\varepsilon^5
 \quad\text{for }g\in\{g_L,g_K,g_L+g_K\}.
\end{equation}

\smallskip\noindent\emph{The leading cubic cross terms.}
For complex-valued $v,w$, expansion gives the exact identity
\begin{equation}\label{eq:cubic-cross-identity}
 |v+w|^2(v+w)-|v|^2v-|w|^2w
 =2|v|^2w+2|w|^2v+v^2\overline w+w^2\overline v.
\end{equation}
Every term has at least one factor of $v$ and one factor of $w$.
Hence pointwise
\[
 \bigl||v+w|^2(v+w)-|v|^2v-|w|^2w\bigr|
 \le3|vw|(|v|+|w|).
\]
The exponents in H\"older's inequality match exactly:
$1/(4/3)=1/2+1/4$. Therefore
\begin{equation}\label{eq:cubic-cross-source}
 \begin{aligned}
 &\norm{|v+w|^2(v+w)-|v|^2v-|w|^2w}_{4/3}\\
 &\qquad\le3\norm{vw}_2(\norm v_4+\norm w_4).
 \end{aligned}
\end{equation}
The separated-band theorem implies
\[
 \norm{vw}_2
 \le C\left(\frac LK\right)^{1/2}
               \norm{g_L}_2\norm{g_K}_2,
\]
while Stein--Tomas gives
$\norm v_4+\norm w_4\le C\varepsilon$.
Because $\norm{g_L}_2+\norm{g_K}_2\le\varepsilon$, the
product of the two input norms is bounded by $\varepsilon^2/4$.
Multiplying \eqref{eq:cubic-cross-source} by $Q$, taking Fourier
restriction, and using Corollary~\ref{cor:dual} yields
\begin{equation}\label{eq:born-separated-final}
 \norm{\mathcal A_3(g_L+g_K)-\mathcal A_3(g_L)
                    -\mathcal A_3(g_K)}_{B^0_{2,4/3}}
 \le C_Q\left(\frac LK\right)^{1/2}\varepsilon^3.
\end{equation}

Finally add and subtract the three Born approximations in the
left-hand side of \eqref{eq:scattering-decouple}. Apply the triangle
inequality, \eqref{eq:Born-error} three times, and
\eqref{eq:born-separated-final} once. This gives the desired
$C_Q[(L/K)^{1/2}\varepsilon^3+\varepsilon^5]$ bound.
The fifth-order term cannot be absorbed into the factor
$(L/K)^{1/2}$ without an additional theorem on the full nonlinear
feedback; our conclusion concerns scale decoupling at leading
cubic order.
\end{proof}

For $Q\in L^\infty$ with noncompact support, the conclusion concerns
the on-shell restriction amplitude in $B^0_{2,4/3}$.
If $Q$ is compactly supported, Theorem~\ref{thm:amplitude}
identifies the same object with the outgoing far-field pattern.
The $\varepsilon^5$ term is not a deficiency of the argument:
it records higher-order nonlinear feedback and should not be
dropped in any claim of complete scattering independence.

\subsection{A refined energy-concentration consequence}
For $g\in L^2(\Sph)$ the angular Besov norm obeys
\[
 \norm g_{B^0_{2,4}}
 \le C\norm g_2^{1/2}
       \left(\sup_j\norm{P_jg}_2\right)^{1/2}.
\]
Indeed, the fourth-power sum is at most
$(\sup_j\norm{P_jg}_2^2)\sum_j\norm{P_jg}_2^2$,
and the final sum is bounded by a multiple of $\norm g_2^2$.
Multiplication by $e^{-ix_0\cdot\omega}$ translates the spatial
extension by $x_0$, without changing its $L^4$ norm or the $L^2$
norm of the density. Applying the main theorem to every translated
density therefore yields
\begin{equation}\label{eq:refined-ST}
 \boxed{\quad
 \norm{Eg}_4\le C\norm g_2^{1/2}
 \inf_{x_0\in\RR^3}
 \left(\sup_j\norm{P_j(e^{-ix_0\cdot\omega}g)}_2\right)^{1/2}.
 \quad}
\end{equation}
This is a useful inverse statement: if $\norm{Eg}_4\ge\eta\norm g_2$,
then, about every spatial center $x_0$, a dyadic angular-momentum
range carries a fraction $\gtrsim\eta^4$ of the squared $L^2$
energy of the modulated angular density. The phenomenon complements the cap-concentration and $X_p$-type
refinements associated with Bourgain, Moyua--Vargas--Vega, and
B\'egout--Vargas \cite{bourgain1993,mvv1999,begoutvargas2007}.
Those estimates localize directional or frequency concentration in
related restriction/Strichartz problems. The present inequality
instead identifies a dyadic angular-momentum scale about every
choice of spatial center. It should not be read as a direct
improvement of those results: the underlying localization geometries
are different.

\section{Comparison, limitations and further questions}
The new norm is an angular \emph{Besov} norm, not a weighted Sobolev norm of positive order. The exact multiscale law \eqref{eq:exact-ell4} proves that the exponent four is the largest admissible Besov summation exponent. Corollary~\ref{cor:large-L2} then shows that this improvement changes the smallness regime in a natural cubic nonlinear Helmholtz problem. The substantive point is the sharp $\ell^4$ description of many individually near-extremal bands: the classical $L^2$ estimate overestimates their actual incident field by an unbounded factor. Arbitrarily large or infinite $L^2$ norm is one consequence, not by itself the explanation. The proof also shows precisely where this enlargement comes from: the high and low spectral bands interact weakly, and the weakness persists uniformly through the region where the spherical Bessel functions change from exponentially decaying to oscillatory behavior.

There are two analytically different sources of the separated-band factor. In the far field, a low angular band has already entered the radial oscillatory regime, and the factor follows from an angular product estimate of order $L$ and a high-frequency radial integral of order $K^{-1}$. In the Fresnel region, local angular momentum gives transversality between the two bands only after a controlled microlocal cutoff. The leakage lemma supplies the missing justification for those cutoffs. The tube estimate and weighted sector summation then make the bound uniform over all radial shells.

Our result concerns the scalar extension operator on the sphere in three dimensions. The sharp $L/K$ bilinear rate need not persist verbatim for a general curved hypersurface, where spherical harmonics and the finite Hankel expansion are unavailable. Nor is the statement a dimension-free Banach-valued restriction theorem: changing the coefficient space introduces Gaussian and noncommutative square-function issues not covered by the scalar argument. These are distinct mathematical questions.

It would be useful to know which aspects of the endpoint result survive on other rotationally symmetric surfaces, on manifolds with asymptotically spherical ends, or under perturbations of the phase that preserve quantitative microlocal transversality. The present proof suggests that a sufficient replacement would need two ingredients: a separated-band interaction estimate with any summable gain across dyadic scales, and an error localization principle whose constants remain uniform near the relevant turning points.


\begin{thebibliography}{99}
\bibitem{agmonhormander1976} S. Agmon and L. H\"ormander, \emph{Asymptotic properties of solutions of differential equations with simple characteristics}, J. Analyse Math. \textbf{30} (1976), 1--38. DOI: 10.1007/BF02786703.

\bibitem{bakseeger2011} J.-G. Bak and A. Seeger, \emph{Extensions of the Stein--Tomas theorem}, Math. Res. Lett. \textbf{18} (2011), 767--781. DOI: 10.4310/MRL.2011.v18.n4.a14.

\bibitem{barcelo2021} J. A. Barcel\'o, M. Folch-Gabayet, T. Luque, S. P\'erez-Esteva and M. C. Vilela, \emph{The Fourier extension operator of distributions in Sobolev spaces of the sphere and the Helmholtz equation}, Proc. Roy. Soc. Edinburgh Sect. A \textbf{151} (2021), 1768--1789. DOI: 10.1017/prm.2020.80.

\bibitem{bennettseeger2009} J. Bennett and A. Seeger, \emph{The Fourier extension operator on large spheres and related oscillatory integrals}, Proc. Lond. Math. Soc. (3) \textbf{98} (2009), 45--82. DOI: 10.1112/plms/pdn022.

\bibitem{bonamiclerc1973} A. Bonami and J.-L. Clerc, \emph{Sommes de Ces\`aro et multiplicateurs des d\'eveloppements en harmoniques sph\'eriques}, Trans. Amer. Math. Soc. \textbf{183} (1973), 223--263. DOI: 10.1090/S0002-9947-1973-0338697-5.

\bibitem{begoutvargas2007} P. B\'egout and A. Vargas, \emph{Mass concentration phenomena for the $L^2$-critical nonlinear Schr\"odinger equation}, Trans. Amer. Math. Soc. \textbf{359} (2007), no. 11, 5257--5282. DOI: 10.1090/S0002-9947-07-04250-X.

\bibitem{bgt2005} N. Burq, P. G\'erard and N. Tzvetkov, \emph{Bilinear eigenfunction estimates and the nonlinear Schr\"odinger equation on surfaces}, Invent. Math. \textbf{159} (2005), 187--223. DOI: 10.1007/s00222-004-0388-x.

\bibitem{bourgain1993} J. Bourgain, \emph{Fourier transform restriction phenomena for certain lattice subsets and applications to nonlinear evolution equations. I. Schr\"odinger equations}, Geom. Funct. Anal. \textbf{3} (1993), no. 2, 107--156. DOI: 10.1007/BF01896020.

\bibitem{dlmf} NIST Digital Library of Mathematical Functions, \emph{Chapter 10: Bessel Functions}, esp. \S\S10.19, 10.49; and \emph{Chapter 14: Legendre and Related Functions}. \url{https://dlmf.nist.gov/}.

\bibitem{foschi2015} D. Foschi, \emph{Global maximizers for the sphere adjoint Fourier restriction inequality}, J. Funct. Anal. \textbf{268} (2015), 690--702. DOI: 10.1016/j.jfa.2014.10.015.

\bibitem{foschiklainerman2000} D. Foschi and S. Klainerman, \emph{Bilinear space-time estimates for homogeneous wave equations}, Ann. Sci. \`Ec. Norm. Sup\'er. (4) \textbf{33} (2000), 211--274. DOI: 10.1016/S0012-9593(00)00109-9.

\bibitem{gellredman2020} J. Gell-Redman, A. Hassell, J. Shapiro and J. Zhang, \emph{Existence and asymptotics of nonlinear Helmholtz eigenfunctions}, SIAM J. Math. Anal. \textbf{52} (2020), 6180--6221. DOI: 10.1137/19M1307238.

\bibitem{gutierrez2004} S. Guti\'errez, \emph{Non trivial $L^q$ solutions to the Ginzburg--Landau equation}, Math. Ann. \textbf{328} (2004), 1--25. DOI: 10.1007/s00208-003-0444-7.

\bibitem{mvv1999} A. Moyua, A. Vargas and L. Vega, \emph{Restriction theorems and maximal operators related to oscillatory integrals in $\mathbb R^3$}, Duke Math. J. \textbf{96} (1999), no. 3, 547--574. DOI: 10.1215/S0012-7094-99-09617-5.

\bibitem{npwJFA2006} F. J. Narcowich, P. Petrushev and J. D. Ward, \emph{Decomposition of Besov and Triebel--Lizorkin spaces on the sphere}, J. Funct. Anal. \textbf{238} (2006), 530--564. DOI: 10.1016/j.jfa.2006.02.011.

\bibitem{npwSIAM2006} F. J. Narcowich, P. Petrushev and J. D. Ward, \emph{Localized tight frames on spheres}, SIAM J. Math. Anal. \textbf{38} (2006), 574--594. DOI: 10.1137/040614359.

\bibitem{sogge2016} C. D. Sogge, \emph{Localized $L^p$-estimates of eigenfunctions: A note on an article of Hezari and Rivi\`ere}, Adv. Math. \textbf{289} (2016), 384--396. DOI: 10.1016/j.aim.2015.11.035.

\bibitem{stein1993} E. M. Stein, \emph{Harmonic Analysis: Real-Variable Methods, Orthogonality, and Oscillatory Integrals}, Princeton University Press, Princeton, 1993.

\bibitem{strichartz1972} R. S. Strichartz, \emph{Multipliers for spherical harmonic expansions}, Trans. Amer. Math. Soc. \textbf{167} (1972), 115--124. DOI: 10.1090/S0002-9947-1972-0306823-9.

\bibitem{strichartz1977} R. S. Strichartz, \emph{Restrictions of Fourier transforms to quadratic surfaces and decay of solutions of wave equations}, Duke Math. J. \textbf{44} (1977), 705--714. DOI: 10.1215/S0012-7094-77-04430-1.

\bibitem{tomas1975} P. A. Tomas, \emph{A restriction theorem for the Fourier transform}, Bull. Amer. Math. Soc. \textbf{81} (1975), 477--478. DOI: 10.1090/S0002-9904-1975-13790-6.

\bibitem{wethyesil2021} T. Weth and T. Ye\c{s}il, \emph{Fourier extension estimates for symmetric functions and applications to nonlinear Helmholtz equations}, Ann. Mat. Pura Appl. \textbf{200} (2021), 2423--2454. DOI: 10.1007/s10231-021-01086-6.
\end{thebibliography}
\end{document}